\RequirePackage{fix-cm}
\documentclass[smallextended]{svjour3}       
\smartqed  
\usepackage{graphicx}
\usepackage{amsmath}
\usepackage{epstopdf} 
\usepackage{subfigure}
\usepackage{amssymb}
\usepackage{amsfonts}
\journalname{Numerical Algorithms}
\begin{document}

\title{Numerical resolution of algebraic equations with symmetries\thanks{This work has been supported by  projects
MTM2010-19510/MTM (MCIN), MTM2012-30860(MECC) and VA118A12-1 (JCYL).
}}


\author{J. Alvarez         \and
        A. Duran 
}


\institute{J. Alvarez \at
              Applied Mathematics Department, University of Valladolid,
Paseo del Cauce 59, E-47011 Valladolid, Spain \\
              \email{joralv@eii.uva.es}           
           \and
           A. Duran \at
              Applied Mathematics Department, University of Valladolid,
Paseo de Bel\'en 15, E-47011 Valladolid, Spain\\
              \email{angel@mac.uva.es}      
}

\date{Received: date / Accepted: date}

\maketitle

\begin{abstract}
Studied here is the effect of the presence of symmetry groups in a system of algebraic equations on the numerical resolution with fixed-point algorithms. It is shown that the symmetries imply two important properties of the system: the solutions are not isolated, but distributed in orbits by the symmetry group and zero is an eigenvalue of the Jacobian evaluated at any of the solutions, being the multiplicity at least the dimension of the symmetry group. The use of the corresponding quotient space by the symmetry group motivates the concept of orbital convergence. This also establishes the conditions under which a fixed-point algorithm converges to some element of the orbit of a solution of the system. Two numerical examples illustrate these results.
\end{abstract}

\section{Introduction}
\label{intro}
The resolution of systems of algebraic equations with symmetry groups appears in the mathematical treatment of many models, \cite{marsdenr,olver}. These groups consist of transformations that map solutions of the algebraic system into other solutions.
Some consequences of the existence of such symmetry groups, from an analytical point of view, have been described in e.~g. \cite{olver}. The main goal of this paper is discussing the influence of the symmetries in the numerical procedures to approximate the solutions of the system. This is surely a question that is more or less implicitly observed but, to our knowledge, not formally discussed, except in particular cases, \cite{alvarezd1}. The study is focused on fixed-point iteration algorithms.

From this numerical setting, two main points of discussion emerge. On the one hand, the formation of the corresponding orbits by the group is an obstacle to assume local uniqueness for the solutions of the system. (This is necessary for some classical local convergence results.) On the other hand, the symmetry group incorporates a zero eigenvalue in the Jacobian of the system at the solution, which implies the existence of the eigenvalue one in the iteration matrix of the corresponding iteration operator if a fixed-point system is considered. This property prevents the application of many results of convergence, \cite{ortegar,ostrowsky,griewank,dennism2}.

More specifically, the following points of the paper are emphasized:
\begin{itemize}
\item Because of the existence of a symmetry group in the system of algebraic equations, the solutions cannot be isolated, since any transformed solution by the symmetry group is also a solution. Therefore, this property forces to talk in terms of orbits of solutions (or orbits of fixed points in the case of a fixed-point system) instead of isolated solutions. Via this symmetry group and the corresponding partition of the space into orbits, the system of equations is associated to a reduced system, whose solutions consist of the orbits of the solutions of the original one.
\item Several points of influence of the symmetry group on the iterative resolution of the system are analyzed. It is first shown that the infinitesimal generators of the group at any solution are eigenvectors of the corresponding Jacobian of the system, associated to the zero eigenvalue. (When the system is of fixed-point type, the infinitesimal generators are eigenvectors associated to the eigenvalue equals one of the corresponding iteration matrix.) When the associated eigenspace is spanned by these infinitesimal generators, the Jacobian can be interpreted in terms of the Jacobian of the reduced system at the corresponding orbit. The rest of the spectrum is preserved by the symmetry group action.

The main consequence, in the numerical resolution, is the concept of orbital convergence. This means local convergence to some element of the orbit of a solution of the system, that is, convergence modulo symmetries of the group. The iterative algorithms are analyzed in the reduced space for the orbit solutions of the reduced system, leading to the derivation of classical convergence results but in the orbital sense. The limit point may be estimated, in a first approximation, by the component of the initial iterate in the eigenspace associated to the eigenvalue equals one of the iteration matrix of the fixed-point iteration function.
\end{itemize}

The structure of the paper is as follows. In Section \ref{sec:1} the definition of a symmetry group of a system of algebraic equations, along with several related concepts, are reminded. The orbit of a point by the group, the infinitesimal generators and the reduced system by the symmetries are emphasized. Some properties of the spectrum of the Jacobian matrix of the system at a solution are described; in particular, the existence of the zero eigenvalue because of the symmetry group and the structure of the associated eigenspace. This information is used in Section \ref{sec:2}, where fixed-point iterations for the numerical resolution of this kind of systems are analyzed. (Now, the zero eigenvalue becomes one for the iteration matrix of the corresponding iteration function. The reduction by symmetries preserves the rest of the spectrum of the iteration matrix.)  Thus, from the introduction of the concept of orbital convergence, classical results, like the Contraction Mapping Theorem, can be stated in an orbital sense. The limit point is an element of the orbit of the fixed point and Section \ref{sec:2} is finished off with some remarks about the structure of the error and estimates on this limit point.  Finally Section \ref{sec:3} is devoted to some numerical illustrations. Two examples are considered. The first system appears in the modelization of an $N$-body problem and involves a one-parameter symmetry group of rotations. The second example concerns the numerical generation of traveling waves in nonlinear dispersive wave equations and takes the Bona-Smith system of Boussinesq equations as a case study.
\section{Systems of algebraic equations with symmetries}
\label{sec:1}
\subsection{Notation and preliminaries}
For details related to the definitions and properties mentioned in this section, see e.~g. \cite{marsdenr,olver}.
Throughout this paper, ${\Omega}$ will stand for a domain in $\mathbb{R}^{m}$ ($m>1$) and $\mathcal{G}$ will denote a $l$-parameter ($l\leq m$) Lie group of transformations
\begin{eqnarray}
\mathcal{G}=\{G_{\alpha}:{\Omega}\rightarrow {\Omega},\alpha=(\alpha_{1},\ldots,\alpha_{l})\in \mathbb{R}^{l}\},\label{na1}
\end{eqnarray}
acting on ${\Omega}$. $\mathcal{G}$ is assumed to satisfy the following conditions:
\begin{itemize}
\item[(A1)] $\mathcal{G}$ is Abelian, connected and global.
\item[(A2)] $\mathcal{G}$ acts regularly on ${\Omega}$.
\end{itemize}
The second property in (A1) means that $\mathcal{G}$ is connected when its structure of manifold is considered and excludes discrete groups in our study. On the other hand, the global action implies that each element $G_{\alpha}\in\mathcal{G}, \alpha\in \mathbb{R}^{l}$ is a diffeomorphism on ${\Omega}$. Finally, it is said that $\mathcal{G}$ acts regularly if the associated orbits
\begin{eqnarray}
\mathcal{G}(x)=\{G_{\alpha}x, \alpha\in \mathbb{R}^{l}\}, \quad x\in \Omega\label{na2}
\end{eqnarray}
are of the same dimension and for each $x\in\Omega$ one can find a neighborhood $W$ such that each orbit intersects $W$ in a pathwise connected subset. This property and the global action imply that the orbits (\ref{na2}) are regular $l$ dimensional submanifolds of ${\Omega}$.

The infinitesimal generators of (\ref{na1}) play a relevant role in the present paper. They are the vector fields
\begin{eqnarray}
g_{j}(x)=\frac{\partial}{\partial \alpha_{j}}G_{\alpha}x\Big|_{\alpha=0},\quad j=1,\ldots,l.\label{na3}
\end{eqnarray}
and provide a local representation of the group, in the sense that the flow associated to each vector field (\ref{na3}) is a one-parameter symmetry group and any element of (\ref{na1}) is locally a composition of these transformations. Note that under the previous hypotheses on $\mathcal{G}$, for each $x\in\Omega$, the vectors $g_{j}(x), j=1,\ldots,l$  are linearly independent and since $\mathcal{G}$ is Abelian, these one-parameter subgroups commute.

Finally the action of (\ref{na1}) can generate the quotient space ${\Omega}/\mathcal{G}$. This consists of identifying points that belong to the same orbit and consequently its elements are orbits in the original manifold ${\Omega}$. The hypotheses assumed on $\mathcal{G}$ imply that the quotient space ${\Omega}/\mathcal{G}$ is a ($m-l$) dimensional manifold.

\subsection{Symmetry groups of equations}
Considered now are systems of algebraic equations on ${\Omega}$ on the form
\begin{eqnarray}
F(x)=0,\quad x\in\Omega,\label{na4}
\end{eqnarray}
where $F=(F_{1},\ldots,F_{m}):{\Omega}\rightarrow\mathbb{R}^{m}$ is assumed to be $C^{1}$ on ${\Omega}$. In the sequel, an equivalent fixed point system
\begin{eqnarray}
x=G(x),\label{na5}
\end{eqnarray}
for some $C^{1}$ iteration function  $G:{\Omega}\rightarrow\Omega$, will also be considered. The equivalence between (\ref{na4}) and (\ref{na5}) is understood 
in the sense that $x^{\ast}\in\Omega$ is a solution of (\ref{na4}) if and only if $x^{\ast}$ is a fixed point of (\ref{na5}), that is,
\begin{eqnarray*}
F(x^{\ast})=0\Leftrightarrow x^{\ast}=G(x^{\ast}).
\end{eqnarray*}
The group (\ref{na1}) is said to be a symmetry group of (\ref{na4}) if each element $G_{\alpha}\in \mathcal{G}$ has the property of transforming solutions of (\ref{na4}) into solutions, that is, \cite{olver}
\begin{eqnarray}
F(x^{\ast})=0\Rightarrow F(G_{\alpha}x^{\ast})=0,\quad \forall \alpha\in \mathbb{R}^{l}.\label{na6}
\end{eqnarray}
Equivalently, the symmetry group condition for (\ref{na5}) is
\begin{eqnarray}
x^{\ast}=G(x^{\ast})\Rightarrow G_{\alpha}x^{\ast}=G(G_{\alpha}x^{\ast}),\quad \forall \alpha\in \mathbb{R}^{l}.\label{na7}
\end{eqnarray}
\subsection{Reduced system}
The presence of the symmetry group (\ref{na1}) has some consequences in the structure of (\ref{na4}) that are relevant for our analysis. Note first that differentiating (\ref{na6}) with respect to $\alpha_{j}, j=1,\ldots,l$ and evaluating at $\alpha=0$ we obtain that if $x^{\ast}$ is a solution of (\ref{na4}) then
\begin{eqnarray*}
F^{\prime}(x^{\ast})g_{j}(x^{\ast})=0,\quad j=1,\ldots,l,\label{na8}
\end{eqnarray*}
where $F^{\prime}(x^{\ast})$ stands for the Jacobian of $F$ at $x^{\ast}$ and the $g_{j}, j=1,\ldots,l$ are given by (\ref{na3}). This means that $F^{\prime}(x^{\ast})$ is singular and the vectors $g_{j}(x^{\ast})$ belong to ${\rm Ker}\left(F^{\prime}(x^{\ast})\right)$. Equivalently, using (\ref{na7}), the iteration matrix $G^{\prime}(x^{\ast})$, corresponding to the fixed point system (\ref{na5}), admits $\lambda=1$ as eigenvalue with the $g_{j}(x^{\ast})$ as eigenvectors.

The second consequence concerns the way the systems (\ref{na4}), (\ref{na5}) can be reduced in the quotient space ${\Omega}/\mathcal{G}$. This can be described by using local coordinates, \cite{olver}. The hypotheses (A1), (A2) allow to construct, in some neighborhood $V$ of $x^{\ast}$, a system of local coordinates $(y,z)=(y_{1},\ldots,y_{l},z_{1},\ldots,z_{m-l})$ for which the group $\mathcal{G}$ corresponds to translations in the first $l$ variables,
\begin{eqnarray*}
G_{\alpha}(y,z)&=&(y_{1}+\alpha_{1},\ldots,y_{l}+\alpha_{l},z_{1},\ldots,z_{m-l}), \\
&&\alpha=(\alpha_{1},\ldots,\alpha_{l})\in \mathbb{R}^{l}, (y,z)\in V.
\end{eqnarray*}
This implies that, locally, the infinitesimal generators are the first $l$ vectors of the canonical basis in $\mathbb{R}^{m}$,
\begin{eqnarray*}
g_{j}(y,z)=\frac{\partial}{\partial \alpha_{j}}G_{\alpha}(y,z)\Big|_{\alpha=0}=e_{j},\quad j=1,\ldots,l,
\end{eqnarray*}
where $e_{j}$ has all the components equals zero except the one in the $j$-th position, which values one. Note then that the symmetry group condition implies that close to $x^{\ast}=(y^{\ast},z^{\ast})$, the system (\ref{na4}) is independent of the first $l$ variables ($F_{i}(y,z)=F_{i}(z), i=1,\ldots,m$) and can be written in the form
\begin{eqnarray}
F_{1}(z_{1},\ldots,z_{m-l})&=&0\nonumber\\
\vdots&&\nonumber\\
F_{l}(z_{1},\ldots,z_{m-l})&=&0\label{na9}\\
F_{1+1}(z_{1},\ldots,z_{m-l})&=&0\nonumber\\
\vdots&&\nonumber\\
F_{m}(z_{1},\ldots,z_{m-l})&=&0,\nonumber
\end{eqnarray}
such that the first $l$ equations in (\ref{na9}) are functionally dependent on the remaining ones. This implies a reduction of the first $l$ variables and the characterization of (\ref{na4}) in terms of the $m-l$ variables $z_{1},\ldots,z_{m-l}$ and the last $m-l$ functionally independent equations, that make up the so-called reduced system. 

In terms of the fixed point system (\ref{na5}) the reduction process means that near the fixed point $x^{\ast}=(y^{\ast},z^{\ast})$ (\ref{na5}) can be written in the form
\begin{eqnarray}
y=y-\widetilde{\bf G}_{1}(z)&&\left\{
\begin{matrix}y_{1}&=&y_{1}-\widetilde{G}_{1}(z_{1},\ldots,z_{m-l})\\
\vdots&&\\
y_{l}&=&y_{l}-\widetilde{G}_{l}(z_{1},\ldots,z_{m-l})
\end{matrix}\right. \label{na10}\\
z=\widetilde{\bf G}_{2}(z)&&\left\{
\begin{matrix}z_{1}&=&\widetilde{G}_{l+1}(z_{1},\ldots,z_{m-l})\\
\vdots&&\\
z_{m-l}&=&\widetilde{G}_{m}(z_{1},\ldots,z_{m-l})
\end{matrix}\right. \label{na11}\
\end{eqnarray}
where (\ref{na11}) is a local representation of the reduced fixed point system and (\ref{na10}) is functionally dependent on (\ref{na11}).

\subsubsection{Example}
In some cases, this reduction can be described more explicitly. This is the case of the following system in $\mathbb{R}^{4}$:
\begin{equation}\label{etna5b}
\begin{array}{l}
F({\bf q})=\omega^{2}{\bf q}+\nabla U({\bf q})=0, {\bf q}=(q_{1},q_{2})^{T}, q_{j}\in \mathbb{R}^{2}, j=1,2,\\
U({\bf q})=-\frac{m_{0}}{|q_{1}|}-\frac{m_{0}}{|q_{2}|}-\frac{1}{|q_{1}-q_{2}|},\;\;
\omega=\sqrt{m_{0}+\frac{1}{4}},
\end{array}
\end{equation}
where $m_{0}>0$. System (\ref{etna5b}) is a particular case of the one considered in the example 1 of section \ref{sec:3}. Some explanations about the related physical model will be then relegated to there, since here we are only interested in illustrating the reduction process. System (\ref{etna5b}) admits the group of rotations
\begin{eqnarray*}
G_{\alpha}{\bf q}=\begin{pmatrix}
\mathcal{G}_{\alpha}q_{1}&0\\0&\mathcal{G}_{\alpha}q_{2}\end{pmatrix},\quad
\mathcal{G}_{\alpha}=\begin{pmatrix}
\cos{\alpha}&-\sin{\alpha}\\
\sin{\alpha}&\cos{\alpha}
\end{pmatrix},\quad \alpha\in\mathbb{R},
\end{eqnarray*}
as a symmetry group. In this case, the reduction can be carried out  by using polar coordinates
\begin{eqnarray*}
q_{1}=(r_{1}\cos\theta_{1},r_{1}\sin\theta_{1}), \quad q_{2}=(r_{2}\cos\theta_{2},r_{2}\sin\theta_{2}).
\end{eqnarray*}
Thus, in terms of the variables
$
{\bf q}\mapsto (r_{1},r_{2},\theta)$, where $\theta=\theta_{1}-\theta_{2}$,
the reduced system is of the form
\begin{eqnarray*}
\omega^{2}r_{1}-\frac{m_{0}}{r_{1}^{2}}-\frac{(r_{1}-r_{2}\cos\theta)}{\left(r_{1}^{2}+r_{2}^{2}-2r_{1}r_{2}\cos\theta\right)^{3/2}}&=&0,\\
\omega^{2}r_{2}-\frac{m_{0}}{r_{2}^{2}}-\frac{(r_{2}-r_{1}\cos\theta)}{\left(r_{1}^{2}+r_{2}^{2}-2r_{1}r_{2}\cos\theta\right)^{3/2}}&=&0,\\
\frac{2r_{1}r_{2}\sin\theta}{\left(r_{1}^{2}+r_{2}^{2}-2r_{1}r_{2}\cos\theta\right)^{3/2}}&=&0.
\end{eqnarray*}

\section{Fixed-point algorithms}
\label{sec:2}
\subsection{Structure of the iteration matrix}
The reduced fixed point system will mainly determine the local behaviour of the corresponding fixed point iteration for (\ref{na5}). Note that due to the local reduced form (\ref{na10}), (\ref{na11}) and the hypotheses on $\mathcal{G}$, in a neighborhood of $x^{\ast}$, the iteration matrix $G^{\prime}(x^{\ast})$ is similar to a block matrix of the form
\begin{eqnarray*}
\begin{pmatrix}
I_{l}&|&-\widetilde{\bf G}_{1}^{\prime}(z^{\ast})\\ 
-&-|-&-\\
0&|&\widetilde{\bf G}_{2}^{\prime}(z^{\ast})
\end{pmatrix}
\end{eqnarray*}
where $x^{\ast}=(y^{\ast},z^{\ast})$, $I_{l}$ stands for the $l\times l$ identity matrix and $\widetilde{\bf G}_{j}^{\prime}(z^{\ast})$ denotes the Jacobian of $\widetilde{\bf G}_{j}(z^{\ast})$ with respect to the variables $z_{1},\ldots,z_{m-l}$ at $z^{\ast}$ for $j=1,2$. Consequently, the spectrum of the iteration matrix $G^{\prime}(x^{\ast})$ consists of the eigenvalue $\lambda=1$, associated to the presence of the symmetry group $\mathcal{G}$ and the spectrum of the iteration matrix of the reduced fixed point system at the orbit of $x^{\ast}$. In particular (see also \cite{durans}), 
\begin{lemma}
\label{lemma21}
Assume the hypotheses (A1), (A2). If
$\lambda=1$ is not an eigenvalue of the iteration matrix of the reduced fixed point system
at the orbit of $x^{\ast}$,
then the geometric and algebraic multiplicities of $\lambda=1$ as eigenvalue of $G^{\prime}(x^{\ast})$ equals $l$ with the vectors in (\ref{na3}) forming a basis.
\end{lemma}
\subsection{Orbital convergence}
The previous results are oriented to analyze the influence of the symmetries in the numerical resolution of (\ref{na4}) with the fixed-point iteration
\begin{eqnarray}
x_{n+1}=G(x_{n}),\quad n=0,1,\ldots,\label{na12}
\end{eqnarray}
and given initial data $x_{0}\in {\Omega}$.

The concept of orbital convergence is now introduced. If $x^{\ast}\in\Omega$ is a solution of (\ref{na4}), the iteration (\ref{na12}) is said to be orbitally convergent to $x^{\ast}$ if there are a neighborhood $B(x^{\ast})$ of $x^{\ast}$ and $\alpha^{\ast}\in \mathbb{R}^{l}$ such that if $x_{0}\in B(x^{\ast})$ then $\{x_{n}\}_{n}$, defined by (\ref{na12}), converges to $G_{\alpha^{\ast}}\left(x^{\ast}\right)$. That is, orbital convergence means that the iteration converges locally to some element of the orbit of $x^{\ast}$ by the symmetry group (\ref{na1}).

Orbital convergence looks to be the natural concept to treat the convergence for this kind of problems, since this must be understood in the quotient manifold ${\Omega}/\mathcal{G}$. The structure of the iteration matrix $G^{\prime}(x^{\ast})$, previously described, clarifies the procedure to obtain results on orbital convergence, since the spectrum of the iteration matrix of the reduced system at the orbit of $x^{\ast}$ coincides with the part of the spectrum of $G^{\prime}(x^{\ast})$ which is not associated to the symmetry group. 
 For instance, the Contraction Mapping Theorem, applied to the reduced system, leads to the following result (see also \cite{alvarezd1}).
\begin{theorem}
\label{theorem1} Assume that Lemma \ref{lemma21} holds and any other eigenvalue of $G^{\prime}(x^{\ast})$, different from $\lambda=1$ is in magnitude below one. Then the iteration (\ref{na12}) converges orbitally to $x^{\ast}$.
\end{theorem}
\subsection{Structure of the error}
The structure of the iteration errors $e_{n}=x_{n}-x^{\ast}, n=0,1,\ldots$ suggests an estimate, at first order, of the element of the orbit of $x^{\ast}$ to which the iteration converges. 
Let us assume that Theorem \ref{theorem1} holds. 
Let $\epsilon>0$ and suppose that the initial error $e_{0}=x_{0}-x^{\ast}$ satisfies $||e_{0}||=O(\epsilon)$. By using the spectrum of $G^{\prime}(x^{\ast})$, we can decompose
\begin{eqnarray*}
e_{0}=\sum_{j=1}^{l}\epsilon\alpha_{j}^{(0)}g_{j}(x^{\ast})+\epsilon v^{(0)},
\end{eqnarray*}
for some $\alpha_{j}^{(0)}\in\mathbb{R},  j=1,\ldots l$. The term $v^{(0)}$ is the component of $e_{0}$ in  the supplementary subspace $V$ of Ker$(I-G^{\prime}(x^{\ast}))$ in $\mathbb{R}^{m}$, consisting of the other $G^{\prime}(x^{\ast})$ invariant generalized eigenspaces associated to the eigenvalues different from $\lambda=1$.
On the other hand, the second term of the error sequence can be written in the form
\begin{eqnarray*}
e_{1}=x_{1}-x^{\ast}=G(x_{0})-G(x^{\ast})=G^{\prime}(x^{\ast})e_{0}+M_{0}
\end{eqnarray*}
where $M_{0}=G(x_{0})-G(x^{\ast})-G^{\prime}(x^{\ast})e_{0}$. If we assume that $||M_{0}||=O(\epsilon^{2})$, \cite{stoerb,dennism1,dennism2,brezinski},
then $e_{1}$ can be written
\begin{eqnarray*}
e_{1}=\sum_{j=1}^{l}\epsilon\alpha_{j}^{(0)}g_{j}(x^{\ast})+\epsilon v^{(1)}+M_{0},
\end{eqnarray*}
where $v^{(1)}=G^{\prime}(x^{\ast})v^{(0)}\in V, ||v^{(1)}||=O(1),  \epsilon\rightarrow 0$. Therefore, the second iterate is of the form
\begin{eqnarray*}
x_{1}=G_{\epsilon\alpha^{(0)}}(x^{\ast})+\epsilon v^{(1)}+\widetilde{M}_{0},
\end{eqnarray*}
with $||\widetilde{M}_{0}||=O(\epsilon^{2})$. The previous arguments can be applied to a general step. We assume that
\begin{eqnarray*}
x_{n}=G_{\epsilon\alpha^{(0)}}(x^{\ast})+\epsilon v^{(n)}+\widetilde{M}_{n-1},
\end{eqnarray*}
with $v^{(n)}\in V, ||v^{(n)}||=O(1),  ||\widetilde{M}_{n-1}||=O(\epsilon^{2}), \epsilon\rightarrow 0$. Then $||e_{n}||=O(\epsilon)$. Now
\begin{eqnarray*}
e_{n+1}=x_{n+1}-x^{\ast}=G(x_{n})-G(x^{\ast})=G^{\prime}(x^{\ast})e_{n}+M_{n}
\end{eqnarray*}
where $M_{n}=G(x_{n})-G(x^{\ast})-G^{\prime}(x^{\ast})e_{n}$. Now, if $||M_{n}||=O(\epsilon^{2})$, then $e_{n+1}$ can be written
\begin{eqnarray*}
e_{n+1}=\sum_{j=1}^{l}\epsilon\alpha_{j}^{(0)}g_{j}(x^{\ast})+\epsilon v^{(n+1)}+M_{n},
\end{eqnarray*}
where $v^{(n+1)}=G^{\prime}(x^{\ast})v^{(n)}\in V, ||v^{(n+1)}||=O(1),  \epsilon\rightarrow 0$. Therefore, the general $(n+1)$-th iterate is 
\begin{eqnarray*}
\label{etna9}
x_{n+1}=G_{\epsilon\alpha^{(0)}}(x^{\ast})+\epsilon v^{(n+1)}+\widetilde{M}_{n},
\end{eqnarray*}
with $||\widetilde{M}_{n}||=O(\epsilon^{2})$. Thus, at $O(\epsilon)$ order, $x_{n}$ approximates the element $G_{\epsilon\alpha^{(0)}}(x^{\ast})$, associated to the projection of the initial iteration $x_{0}$  onto Ker$(I-G^{\prime}(x^{\ast}))$, plus a term $v^{(n)}$ that, according to Theorem \ref{theorem1}, converges to zero.

\section{Numerical experiments}
\label{sec:3}
The previous results are illustrated by two examples.
\subsection{Example 1. A $N$-body problem}
\begin{figure}[htbp]
\centering
\subfigure[]{
\includegraphics[width=0.45\textwidth]{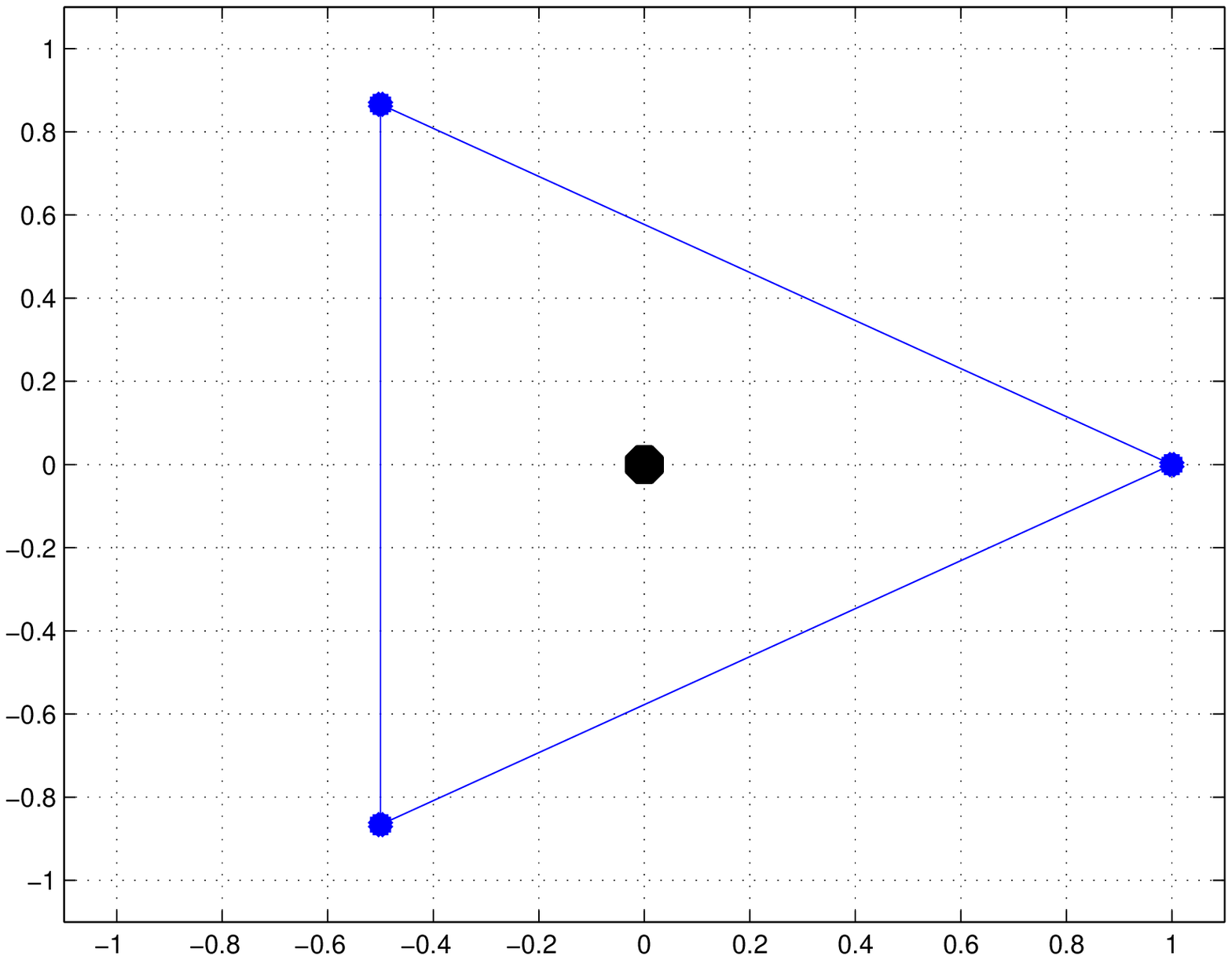} }
\subfigure[]{
\includegraphics[width=0.45\textwidth]{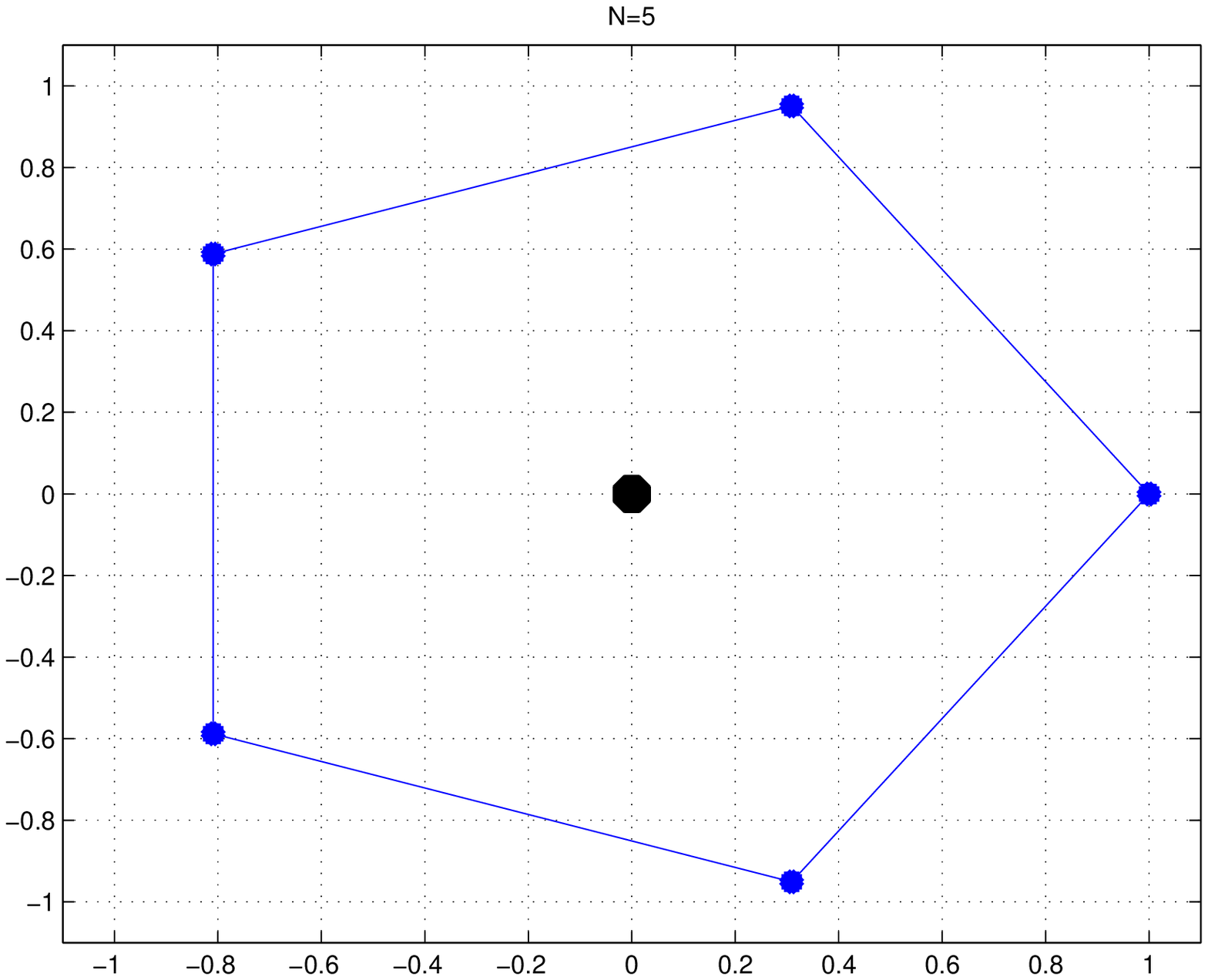} }
\subfigure[]{
\includegraphics[width=0.45\textwidth]{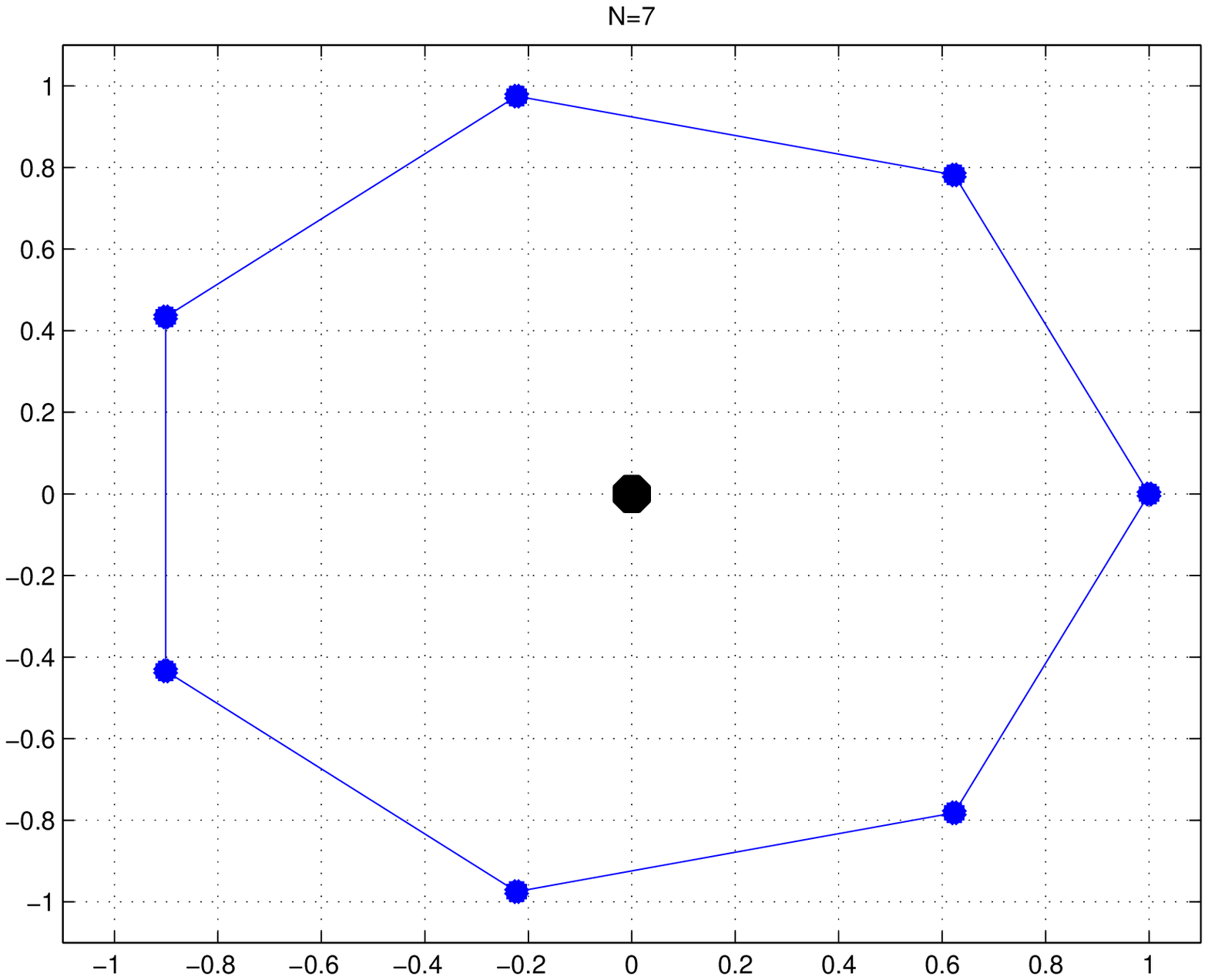} }
\caption{Planar configurations (\ref{etna42}): (a) $N=3$, (b) $N=5$, (c) $N=7$.\label{nbody1}}
\end{figure}

Given $N>1$, $m_{j}>0, j=0,\ldots, N$ and
\begin{eqnarray*}
\omega=\sqrt{m_{0}+\frac{1}{4}\sum_{k=1}^{N-1}{\rm csc}\frac{\pi k}{N}},
\end{eqnarray*}
we consider the $2N\times 2N$ system
\begin{eqnarray}
\omega^{2}M {\bf q}=-\nabla U({\bf q}),\label{etna41}
\end{eqnarray}
where $M=diag(m_{1},m_{1},\ldots,m_{N},m_{N})$ is a $2N\times 2N$ diagonal matrix with the indicated diagonal elements, ${\bf q}=(q_{1}^{T},\ldots,q_{N}^{T})^{T}$ with $q_{j}\in \mathbb{R}^{2}, j=1,\ldots N$ and
\begin{eqnarray}
U({\bf q})=-\sum_{j=1}^{N}\frac{m_{0}m_{j}}{|q_{j}|}-\sum_{i<j}^{N}\frac{m_{i}m_{j}}{|q_{i}-q_{j}|}.\label{etna41b}
\end{eqnarray}
System (\ref{etna41}), (\ref{etna41b}) appears in the modelization of the $N$-body problem, where $q_{j}$ stands for the position of the $j$-th body with mass $m_{j}$ (assuming a planar configuration) and where $m_{0}$ denotes a central mass, which is supposed to be at the origin, \cite{arnold}. (The example considered in section \ref{sec:1} is a particular case with $N=2, m_{1}=m_{2}=1$. This will be used again in the numerical experiments below.) In the case of (\ref{etna41}), the solutions correspond to a configuration where bodies form regular polygons around the mass $m_{0}$, see Figure \ref{nbody1}. Assuming for simplicity that $m_{j}=1, j=1,\ldots N$, a solution of (\ref{etna41}), (\ref{etna41b}) is of the form, \cite{moeckel}
\begin{eqnarray}
{\bf q}^{*}=(q_{1}^{*},\ldots,q_{N}^{*})^{T},\quad q_{j}^{*}=(\cos\theta_{j},\sin\theta_{j})^{T}, \theta_{j}=\frac{2\pi j}{N}, \quad j=1,\ldots,N.\label{etna42}
\end{eqnarray}
The dynamics of the corresponding solution of the $N$-body system with initial data given by ${\bf q}^{*}$ corresponds to a uniform rotation around $m$ with angular velocity given by $\omega$, \cite{canod}.

\begin{figure}[htbp]
\centering
\includegraphics[width=0.6\textwidth]{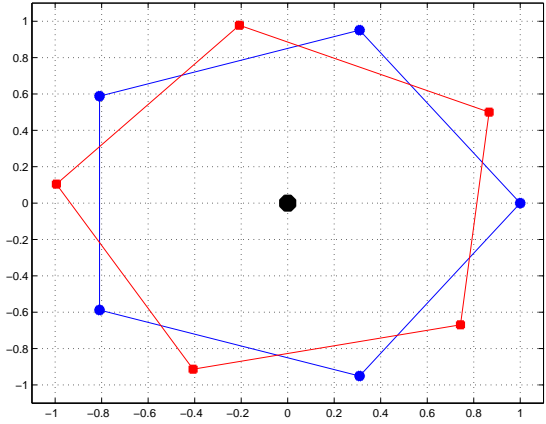}
\caption{Rotation (\ref{etna44}) with $\alpha=\pi/6$. \label{nbody2}}
\end{figure}

System (\ref{etna41}), (\ref{etna41b}) admits the symmetry group of rotations (see Figure \ref{nbody2})
\begin{eqnarray}
\mathcal{G}=\{G_{\alpha}:\mathbb{R}^{2N}\rightarrow\mathbb{R}^{2N},\alpha\in\mathbb{R}\},\label{etna43}
\end{eqnarray}
where, for ${\bf q}\in \mathbb{R}^{2N}, {\bf q}=(q_{1}^{T},\ldots,q_{N}^{T})^{T}$, $q_{j}\in\mathbb{R}^{2}, j=1,\ldots,N$,
\begin{eqnarray}
\label{etna44}
G_{\alpha}{\bf q}=diag\left(
\mathcal{G}_{\alpha}q_{1},\ldots,\mathcal{G}_{\alpha}q_{N}\right),\quad
\mathcal{G}_{\alpha}=\begin{pmatrix}
\cos{\alpha}&-\sin{\alpha}\\
\sin{\alpha}&\cos{\alpha}
\end{pmatrix}.
\end{eqnarray}
The infinitesimal generator of the group is
\begin{eqnarray}
\label{etna44b}
v({\bf q})=\frac{d}{d\alpha}G_{\alpha}{\bf q}\big|_{\alpha=0}=diag\left(Jq_{1},\ldots,Jq_{N}\right),\quad J=\begin{pmatrix}
0&-1\\1&0\end{pmatrix}.
\end{eqnarray}
By way of first illustration of the orbital convergence, an iterative method for (\ref{etna41}), (\ref{etna41b}) is implemented with the goal of approximating the exact solution (\ref{etna42}) and its corresponding orbit. Note that system (\ref{etna41}) can be written in the form (\ref{na5}) with
\begin{eqnarray}
G({\bf q})=-\frac{1}{\omega^{2}}M^{-1}\nabla U({\bf q}).\label{etna44c}
\end{eqnarray}
Observe that $G$ is an homogeneous function of ${\bf q}$ with degree $p=-2$. In particular, this implies that
\begin{eqnarray*}
G^{\prime}({\bf q}^{*}){\bf q}^{*}=pG({\bf q}^{*})=p{\bf q}^{*},
\end{eqnarray*}
and therefore $p=-2$ is an eigenvalue of the iteration matrix $G^{\prime}({\bf q}^{*})$. This means that the fixed-point algorithm with iteration function $G$ will not be convergent in general. A typical alternative for this kind of homogeneous systems is the so-called Petviashvili method. This method was originated in the context of the generation of lump solitary waves in the KP-I equation, \cite{petviashvili}, and it is widely used in this context, where nonlinearities with homogeneous functions are frequent. The application of the method to (\ref{etna41}), (\ref{etna41b}) is as follows. From ${\bf q}_{0}\neq 0$, a sequence of approximations to ${\bf q}^{*}$ is generated by the formulas
\begin{eqnarray}
s_{n}&=&\frac{\langle\omega^{2} M{\bf q}_{n},{\bf q}_{n}\rangle}{\langle -\nabla U({\bf q}_{n}),{\bf q}_{n}\rangle},\label{etna45}\\
{\bf q}_{n+1}&=&G_{Pet}({\bf q}_{n})=-\frac{s_{n}^{\gamma}}{\omega^{2}}M^{-1}\nabla U({\bf q}_{n}),\quad n=0,1,\ldots,\label{etna46}
\end{eqnarray}
where $\gamma$ is a free parameter and $\langle\cdot,\cdot\rangle$ stands for the Euclidean inner product in $\mathbb{R}^{2N}$.  The analysis of the method has been carried out in several papers and for some cases \cite{pelinovskys,lakobay}. Essentially, local convergence, in both the classical and orbital sense, is obtained when $\lambda=p$ is simple and is the only eigenvalue of $G^{\prime}({\bf q}^{*})$ with modulus greater than one and $\gamma$ is chosen under the condition $|p+\gamma (1-p)|<1$, with the best rate of convergence when $\gamma=p/(p-1)$. The arguments to obtain this result rely on the effect of the so-called stabilizing factor (\ref{etna45}) on the spectrum of $G^{\prime}({\bf q}^{*})$. This factor acts like a filter, eliminating the direction of divergence provided by the eigenspace associated to $\lambda=p$, \cite{alvarezd1}. Thus, the Jacobian of the iteration function of the procedure (\ref{etna45}), (\ref{etna46}) at the fixed point shares the spectrum of $G^{\prime}({\bf q}^{*})$ except the \lq harmful\rq\ (from the point of view of the convergence) eigenvalue $\lambda=p$, which is transformed to an eigenvalue with modulus below one (eventually zero in the case $\gamma=p/(p-1)$).

The procedure (\ref{etna45}), (\ref{etna46}) is written as a fixed-point algorithm with iteration function given by $G_{Pet}$. It will be considered to illustrate the orbital convergence for the system (\ref{etna41}), (\ref{etna41b}). First, and according to the previous comments, the convergence can be checked from Tables \ref{tab_na4} and \ref{tab_na1}. For $N=2$, $m_{j}=1, j=1,2$, and several values of $m_{0}$, these tables  display, respectively, the eigenvalues
of $G^{\prime}({\bf q}^{*})$ and of the iteration matrix  $G_{Pet}^{\prime}({\bf q}^{*})$, evaluated at the exact ${\bf q}^{*}$ given by (\ref{etna42}).
\begin{table}
\begin{center}
\begin{tabular}{|c|c|c|c|c|}
\hline\hline
$m_{0}=10$ & $m_{0}=5$&$m_{0}=4$&$m_{0}=1$&$m_{0}=0$\\\hline
-2.0000E+00&-2.0000E+00&-2.0000E+00&-2.0000E+00&-2.0000E+00\\
9.9999E-01&1.0000E+00&1.0000E+00&-1.6000E+00&-2.0000E+00\\
-5.7143E-01&-8.8888E-01&-1.0000E+00&1.0000+00&1.0000E+00\\
2.8571E-01&4.4444E-01&5.0001E-01&8.0000E-01&1.0000E+00\\
\hline\hline
\end{tabular}
\end{center}
\caption{Eigenvalues of the iteration matrix
$G^{\prime}({\bf q}^{*})$ for several values of $m_{0}$.}\label{tab_na4}
\end{table}
\begin{table}
\begin{center}
\begin{tabular}{|c|c|c|c|c|}
\hline\hline
$m_{0}=10$ & $m_{0}=5$&$m_{0}=4$&$m_{0}=1$&$m_{0}=0$\\\hline
9.9999E-01&1.0000E+00&1.0000E+00&-1.6000E+00&-2.0000E+00\\
-5.7143E-01&-8.8888E-01&-1.0000E+00&1.0000+00&1.0000E+00\\
2.8571E-01&4.4444E-01&5.0000E-01&8.0000E-01&1.0000E+00\\
2.2205E-16&5.5511E-17&-5.5511E-16&0&0\\
\hline\hline
\end{tabular}
\end{center}
\caption{Eigenvalues of the iteration matrix
$G_{Pet}^{\prime}({\bf q}^{*})$ for several values of $m_{0}$.}\label{tab_na1}
\end{table}

In the case of Table \ref{tab_na4}, the degree of homogeneity $p=-2$ appears as simple, dominant eigenvalue, as well as $\lambda=1$, also simple. This last one corresponds to the symmetry group of rotations (\ref{etna43}), (\ref{etna44}). The effect of the method (\ref{etna45}), (\ref{etna46}) on the spectrum of the iteration matrix is observed in Table \ref{tab_na1}. The eigenvalue $p=-2$ is filtered and becomes zero in the spectrum of the new matrix, while the rest of the eigenvalues is preserved.

\begin{figure}[htbp]
\centering \subfigure[]{
\includegraphics[width=0.45\textwidth]{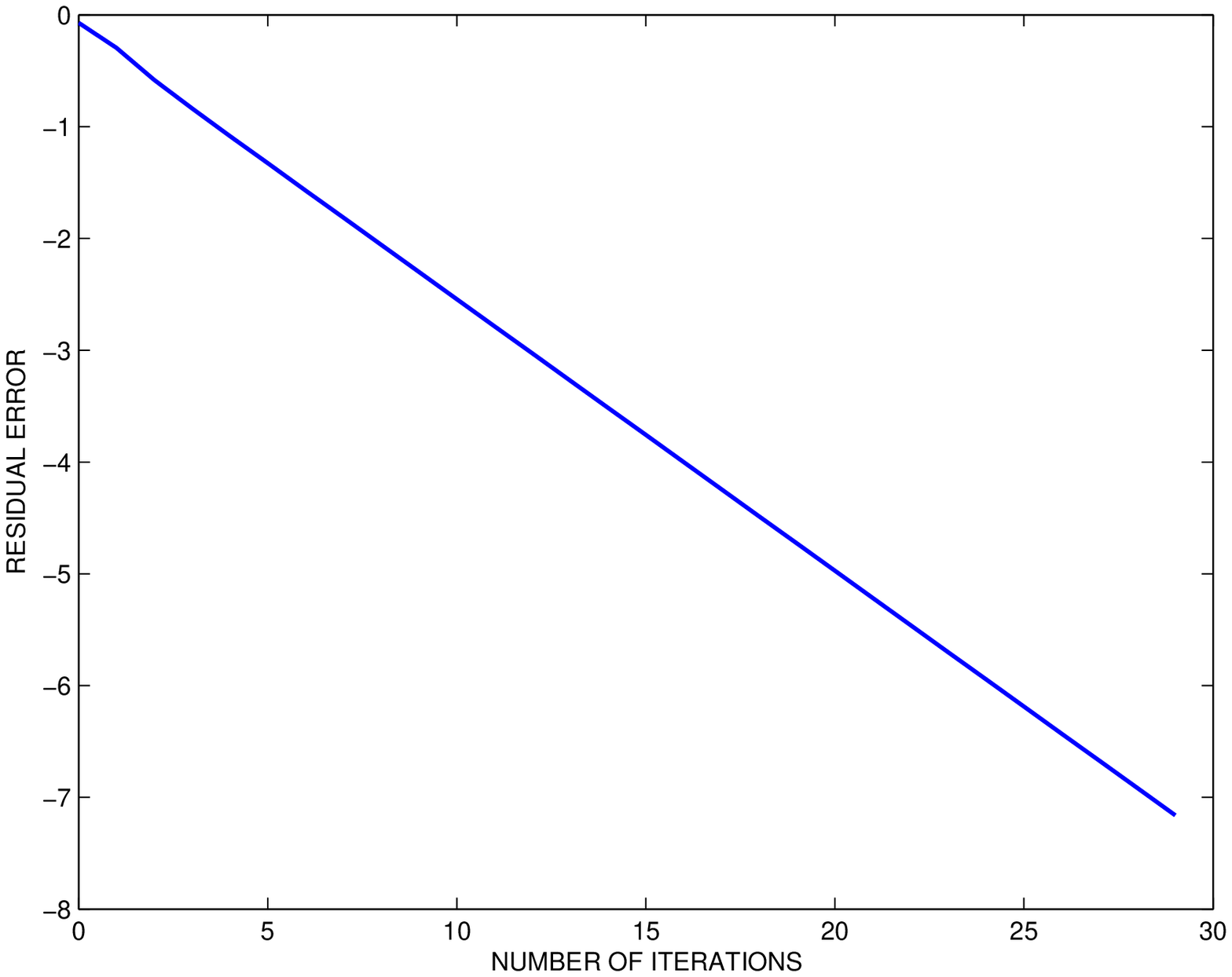} }
\subfigure[]{
\includegraphics[width=0.45\textwidth]{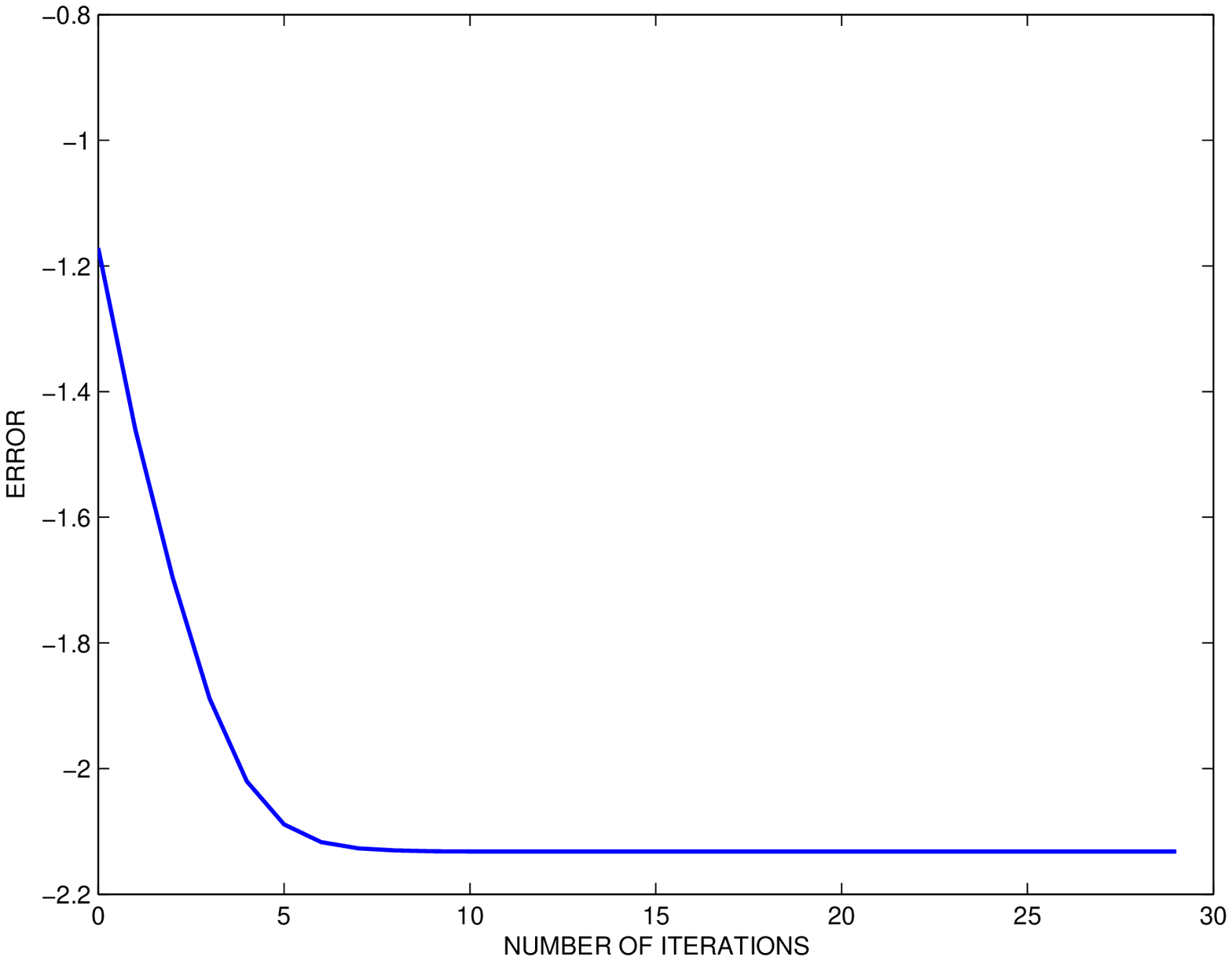} }
\subfigure[]{
\includegraphics[width=0.45\textwidth]{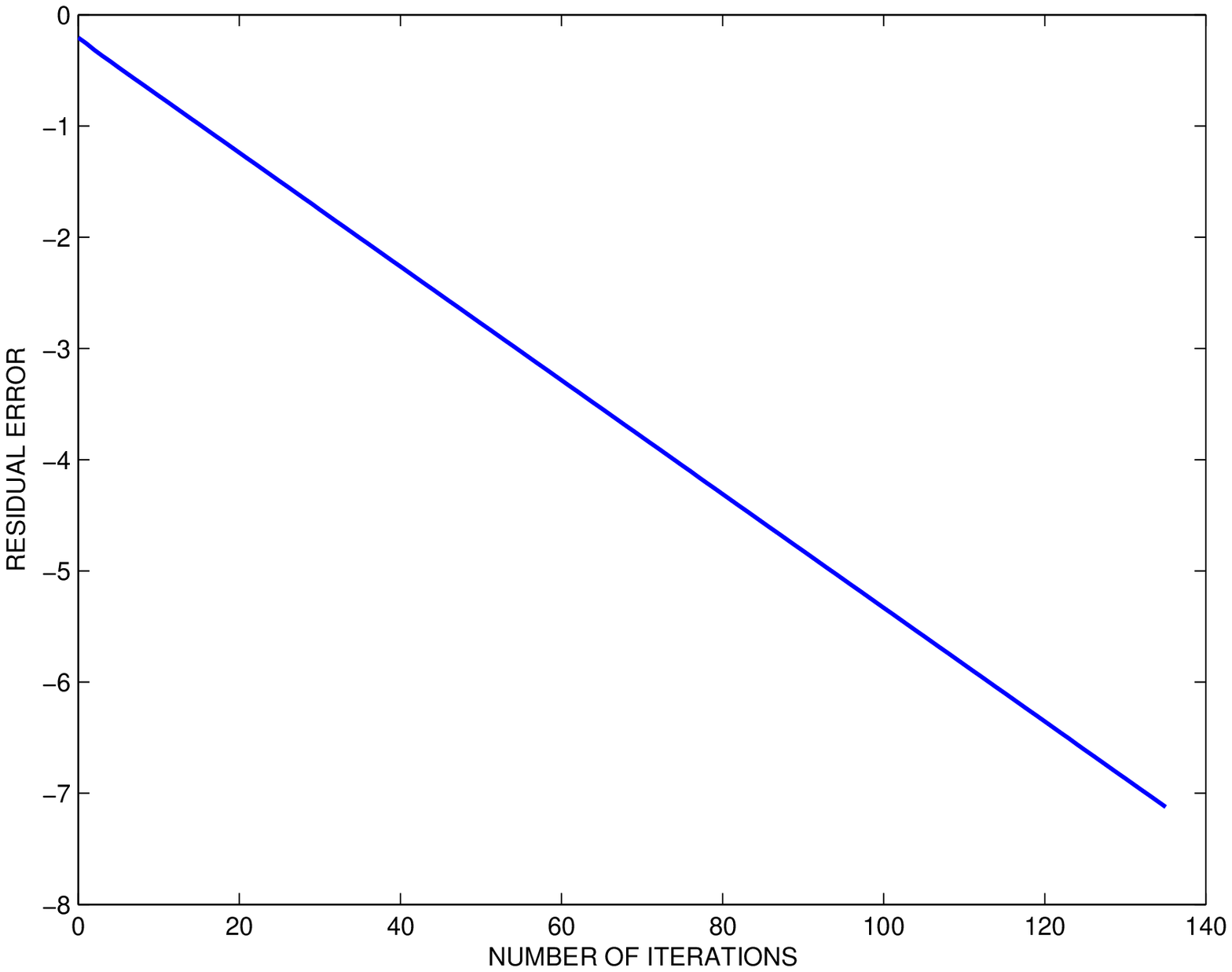} }
\subfigure[]{
\includegraphics[width=0.45\textwidth]{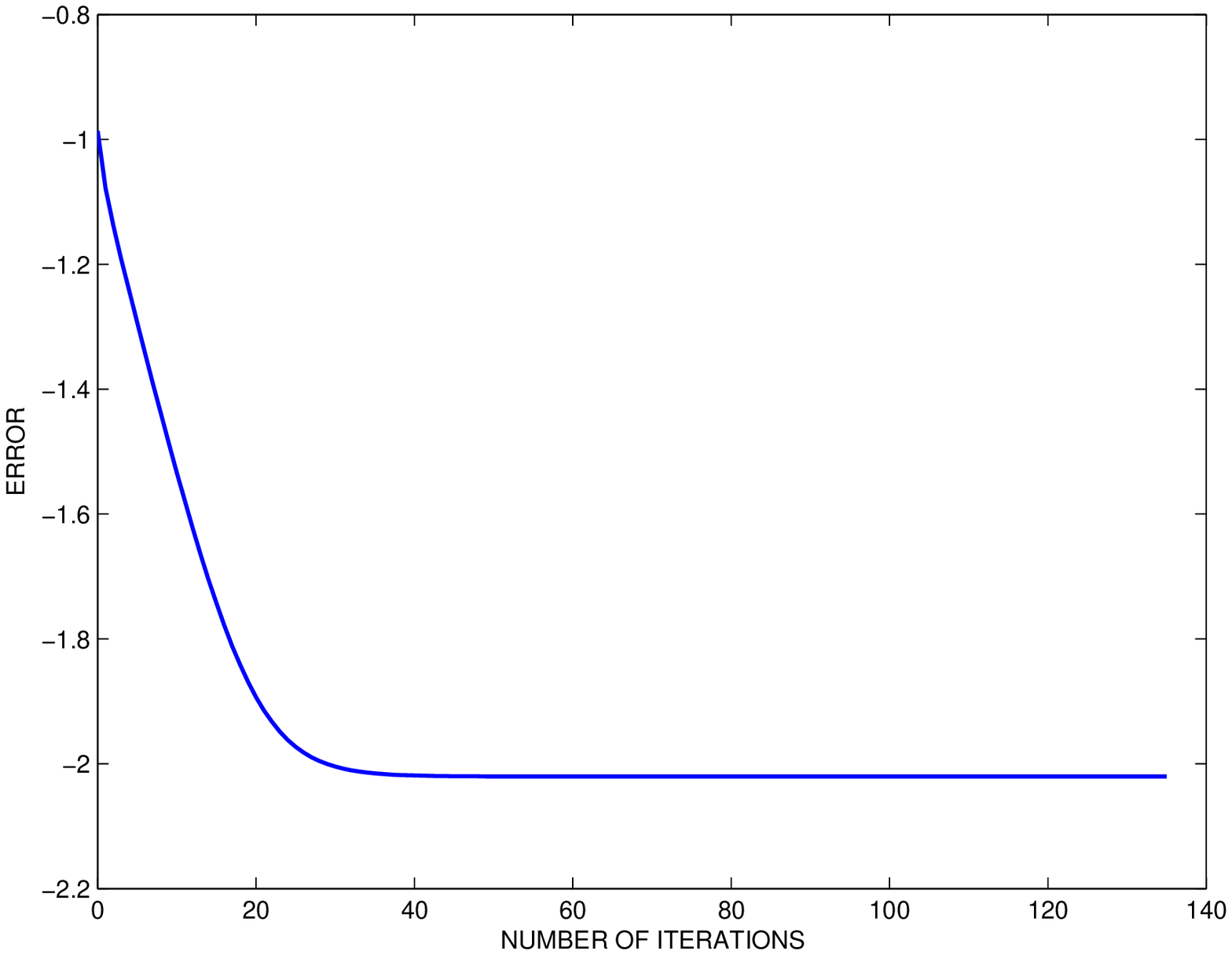} }
\caption{Convergent case: residual error and error with respect to the exact $q^{*}$ as functions of number of iterations: (a), (b) $m_{0}=10$; (c), (d) $m_{0}=5$. The initial iteration is given by (\ref{etna46b}) with $\epsilon=0.1$.} \label{rtb1}
\end{figure}

Convergence is now discussed. According to the parameter $m_{0}$, the following results include two convergent cases
($m_{0}=10, 5$), one not convergent case ($m_{0}=4$) and two
divergent cases ($m_{0}=1, 0$). The convergent case is
illustrated in Figure \ref{rtb1} , where the initial iteration
is a perturbation of ${\bf q}^{*}$ of the form
\begin{eqnarray}
{\bf q}_{0}={\bf q}^{*}+\epsilon w,\label{etna46b}
\end{eqnarray}
with $\epsilon=0.1$ and $w$ is the vector of ones of length $2N$. The procedure (\ref{etna45}), (\ref{etna46}) is run until one of the following errors is below a tolerance $TOL=1E-07$:
\begin{itemize}
\item[(i)] The residual error
\begin{eqnarray}
RE_{n}=||{\bf q}_{n}-G({\bf q}_{n})||,\quad n=0,1,\ldots,\label{RE}
\end{eqnarray} where $G$ is given by (\ref{etna44c}).
\item[(ii)] The (relative) error with respect to the exact fixed point (\ref{etna42})
\begin{eqnarray}
E_{n}=||{\bf q}_{n}-{\bf q}^{\ast}||/||{\bf q}^{\ast}||.\label{ER}
\end{eqnarray}
\end{itemize}
In both cases, the Euclidean norm is used and the errors are displayed in logarithmic scale and as functions of the number of iterations. Figures \ref{rtb1} (a) and (b) show (\ref{RE}) and (\ref{ER}), respectively, for the case $m_{0}=10$ while Figures \ref{rtb1} (c) and (d) correspond to $m_{0}=5$. We observe that the larger $m_{0}$ the faster the iteration converges. The behaviour of the residual error (\ref{RE}) shows the convergence of the iteration. On the other hand, the errors (\ref{ER}), in both cases, stop decreasing from a certain iteration. This may indicate that this convergence is probably to an element of the orbit of ${\bf q}^{\ast}$, different from it and generated by a small component, in the initial iteration, associated to the eigenspace of the eigenvalue $\lambda=1$.

\begin{figure}[htbp]
\centering
\includegraphics[width=8cm]{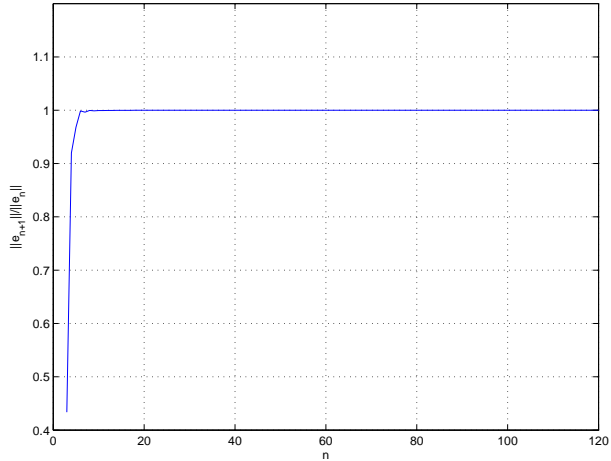}
\caption{Ratios $E_{n+1}/E_{n}$ VS number of iterations ($m_{0}=5$) for (\ref{etna46}) with ${\bf q}_{0}$ given by (\ref{etna46b}).} \label{rtb2b}
\end{figure}

The linear order of convergence is confirmed in Figure \ref{rtb2b}, where the ratios $E_{n+1}/E_{n}$ are displayed as function of the number $n$ of iterations. As $n$ grows, the quotients tend to one.

\begin{figure}[htbp]
\centering \subfigure[]{
\includegraphics[width=0.7\textwidth]{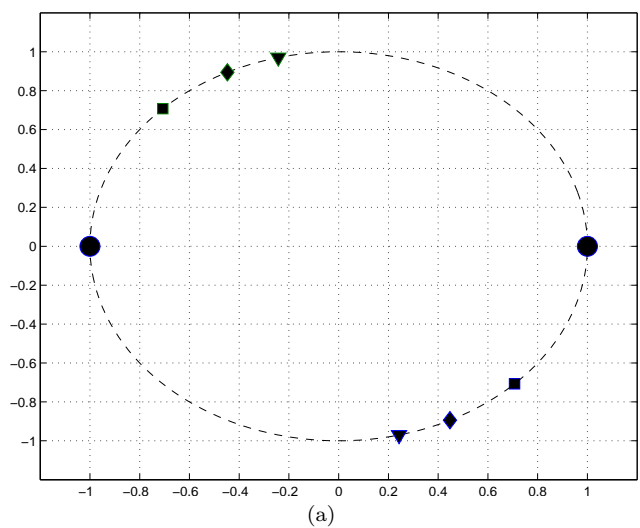} }
\subfigure[]{
\includegraphics[width=0.7\textwidth]{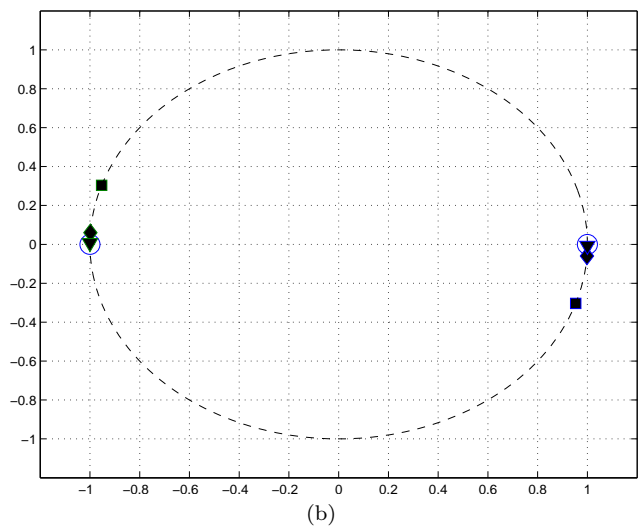} }
\caption{Iteration of system (\ref{etna41}). (a) Final iterates from (\ref{etna47}) with $\epsilon=1$ (with squares), $\epsilon=2$ (diamonds) and $\epsilon=4$ (triangles). (b) Final iterates from (\ref{etna46b}) with $\epsilon=1$ (with squares), $\epsilon=2$ (diamonds) and $\epsilon=4$ (triangles).} \label{rtb3}
\end{figure}

The orbital convergence and the dependence on the initial iteration are also illustrated in Figures \ref{rtb3}(a) and (b). Figure \ref{rtb3}(a) displays the position of the fixed point ${\bf q}^{*}$ (with filled circles) and of three final iterates of the algorithm (\ref{etna45}), (\ref{etna46}) obtained from initial iterations of the form
\begin{eqnarray}
{\bf q}_{0}={\bf q}^{*}+\epsilon v({\bf q}^{*}),\label{etna47}
\end{eqnarray}
where $v({\bf q}^{*})$ is given by (\ref{etna44b}) and for $\epsilon=1$ (with squares), $\epsilon=2$ (diamonds) and $\epsilon=4$ (triangles). The fixed mass is $m_{0}=10$.  Since $v$ is the infinitesimal generator of the symmetry group (\ref{etna43}), the different choices (\ref{etna47}) reflect the influence of ${\bf q}_{0}$ on the limit of the iteration, which in the three cases belongs to the orbit of ${\bf q}^{*}$ by the symmetry group. This orbital convergence is also observed in Figure \ref{rtb3}(b), which shows the final iterates obtained from initial iterations of the form (\ref{etna46b}) where $\epsilon=1$ (with squares), $\epsilon=2$ (diamonds) and $\epsilon=4$ (triangles).
\begin{figure}[htbp]
\centering \subfigure[]{
\includegraphics[width=0.45\textwidth]{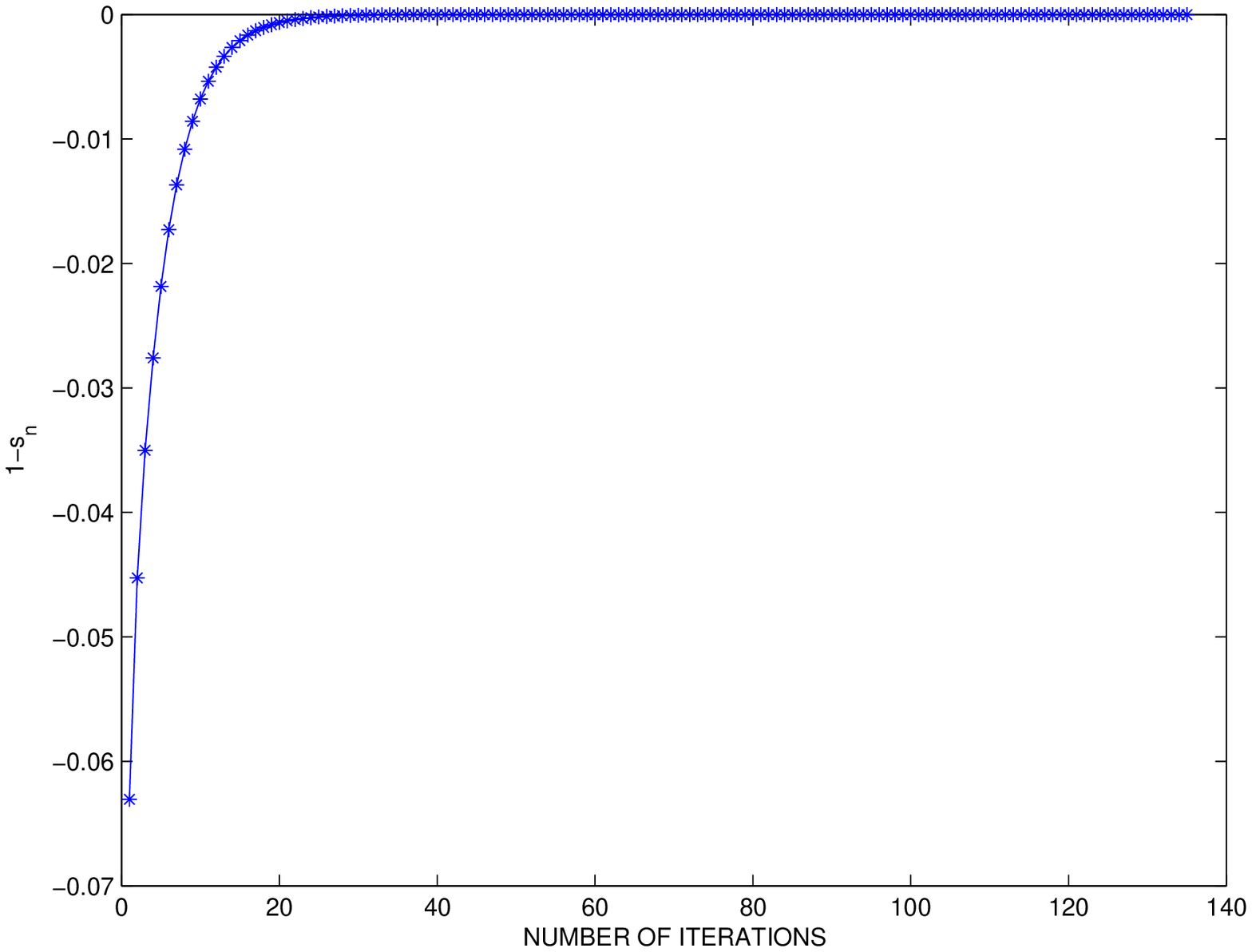} }
\subfigure[]{
\includegraphics[width=0.45\textwidth]{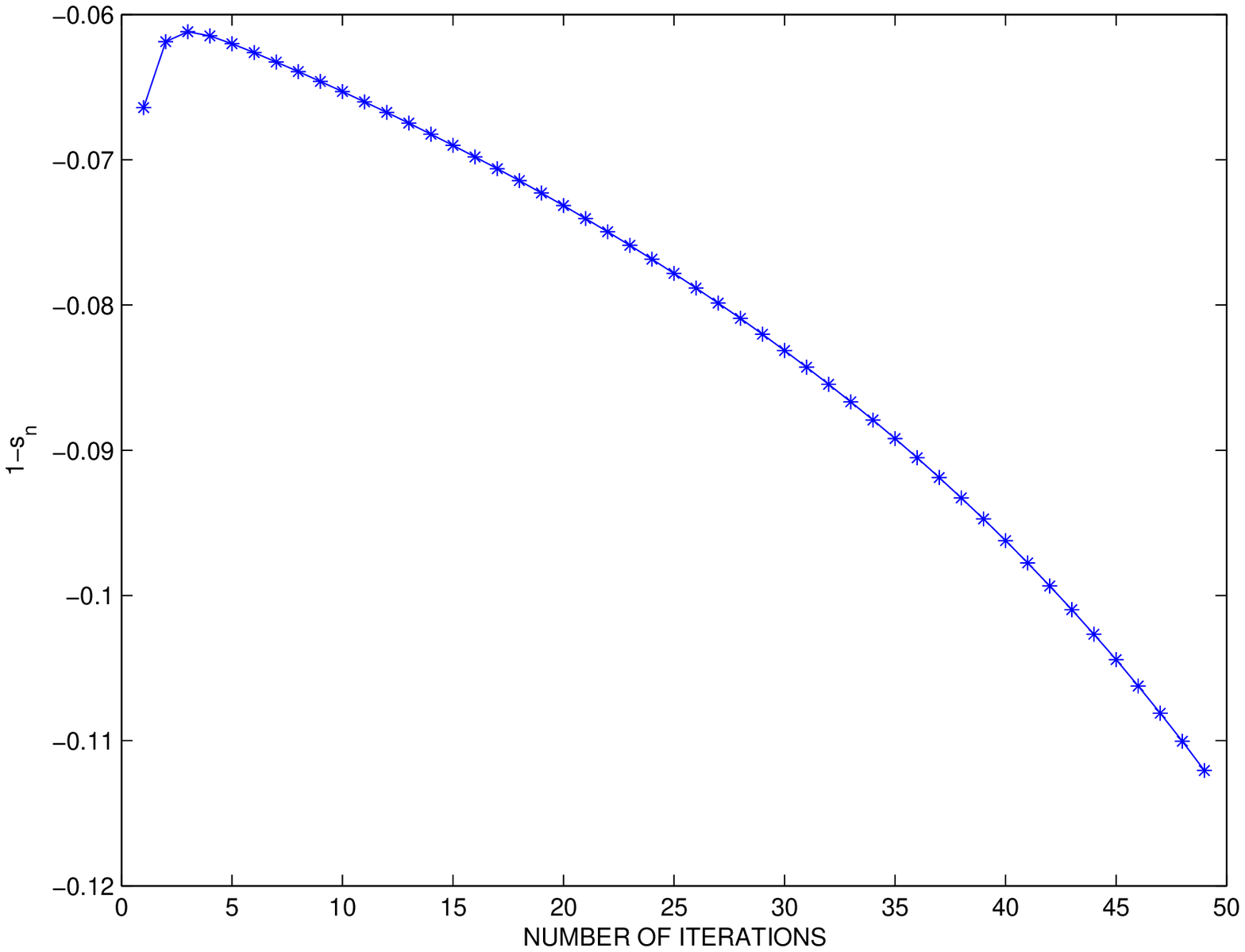} }
\subfigure[]{
\includegraphics[width=0.45\textwidth]{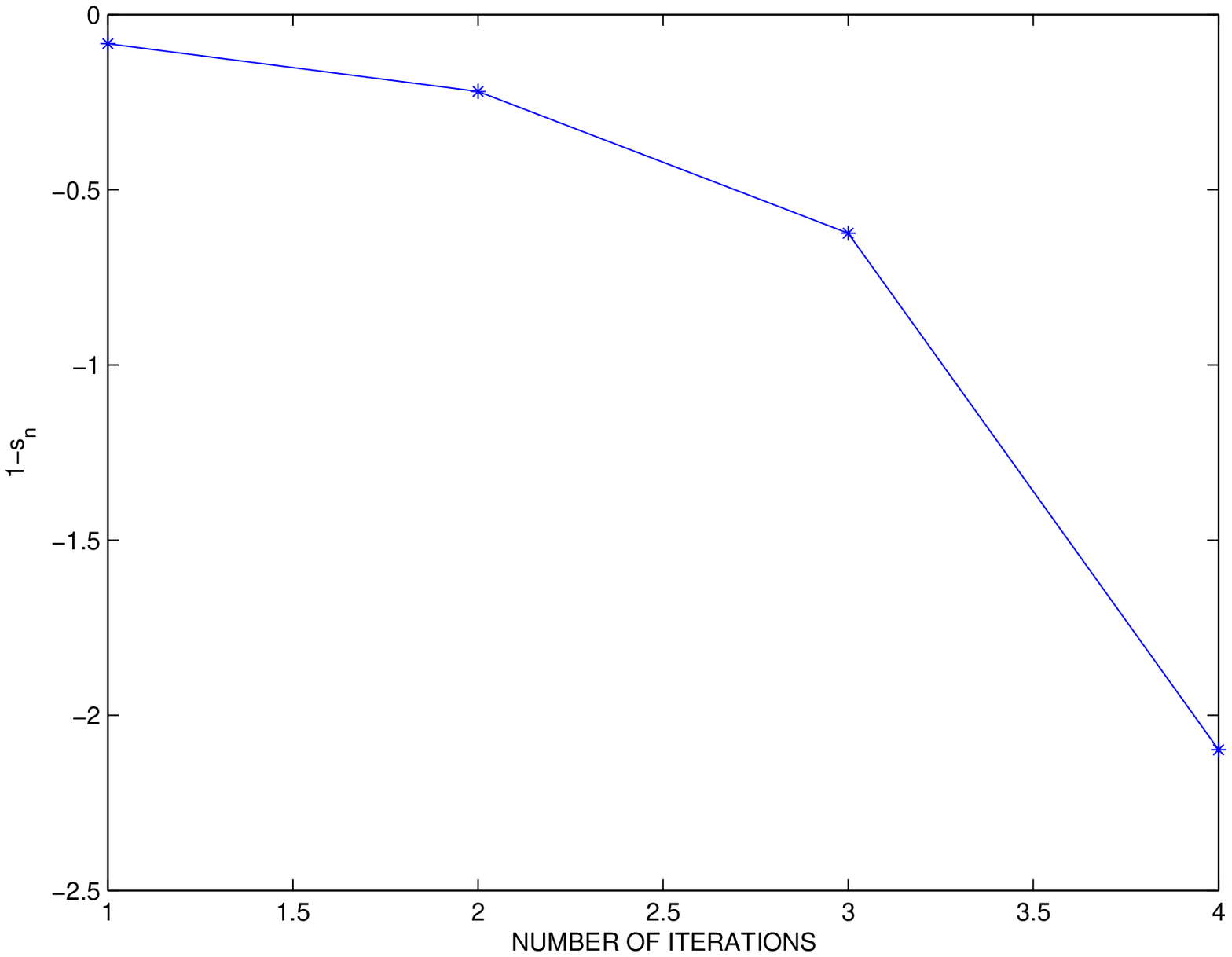} }
\subfigure[]{
\includegraphics[width=0.45\textwidth]{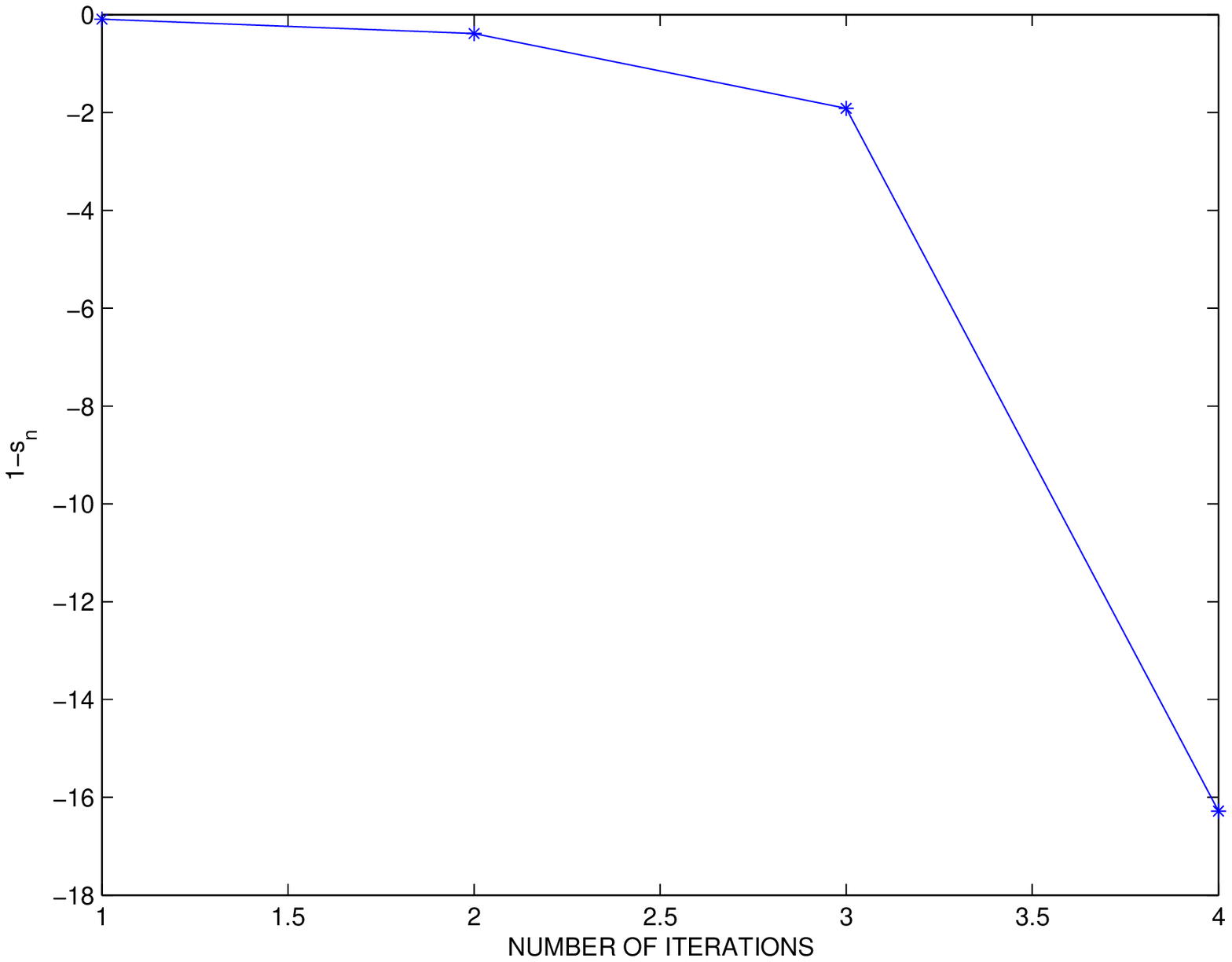} }
\caption{Discrepancy $1-s_{n}$ of the stabilizing factor vs number of iterations: (a) $m_{0}=5$ (convergent);  (b) $m_{0}=4$ (not convergent); (c) $m_{0}=1$ (divergent); (d) $m_{0}=0$ (divergent). The initial iteration is given by (\ref{etna46b}) with $\epsilon=0.1$.} \label{rtb2}
\end{figure}

The behaviour of the stabilizing factor (\ref{etna45}), for several situations, is illustrated in Figure \ref{rtb2}. Note that (\ref{etna45}), when evaluated at the fixed point ${\bf q}^{\ast}$, is equal to one. Thus, in the case of convergence, the factors $s_{n}$ must tend to one. This is observed in Figure \ref{rtb2}(a), that shows the discrepancy $1-s_{n}$ as function of the number of iterations and for the case $m_{0}=5$. A different behaviour of the sequence of stabilizing factors indicates nonconvergence of (\ref{etna46}). For example, the behaviour in Figure \ref{rtb2}(b), corresponding to $m_{0}=4$, may be justified by the presence, in the initial iteration (\ref{etna46b}), of a relevant component in the direction of the eigenspace associated to the eigenvalue $-1$ (see Tables \ref{tab_na4} and \ref{tab_na1}, third column). Note that the eigenvalue $\lambda=1$ does not affect the behaviour of the stabilizing factor, since the quotient in (\ref{etna45}) equals one at any element of the orbit of ${\bf q}^{\ast}$.
Finally, the divergence in Figures \ref{rtb2} (c) and (d) occurs
because of different reasons, see Tables \ref{tab_na4} and \ref{tab_na1}, fourth and fifth columns. In the first case, there is an
additional eigenvalue of modulus greater than one and different
from $p=-2$; in the second case, the dominant eigenvalue is not
simple. This last situation corresponds somehow to a different model. When $m_{0}=0$, the system is affected by forces which are inversely proportional to the distance between the bodies, see (\ref{etna41b}). Its inclusion, however, illustrates a different divergent situation. Note also the presence of a double eigenvalue $\lambda=1$.

\subsection{Example 2. Traveling wave generation}
The second example involves the translational symmetry with the numerical generation of solitary-wave solutions of the Bona-Smith system \cite{bonacs1,bonacs2,boussinesq,bonas}
\begin{equation}\label{bse}
\begin{array}{l}
\eta_{t}+u_{x}+(\eta u)_{x}-b\eta_{xxt}=0,\\
u_{t}+\eta_{x}+uu_{x}+c\eta_{xxx}-du_{xxt}=0,
\end{array}
\end{equation}
for $x\in \mathbb{R},
t\geq 0$, with initial conditions
$$
\eta(x,0)=\eta_{0}(x),\quad u(x,0)=u_{0}(x), \quad x\in\mathbb{R}
$$ and
$$
c=\frac{2}{3}-\theta^{2},
b=d=\frac{1}{2}(\theta^{2}-\frac{1}{3}),\theta^{2}\in
(\frac{2}{3},1].
$$
System (\ref{bse}) is a particular case of the family of Boussinesq systems considered in water wave theory, introduced in \cite{bonacs1,bonacs2}. It is an approximation to the two-dimensional Euler equations for the irrotational free surface flow of an incompressible, inviscid fluid in a uniform, horizontal channel, when no cross channel variations are considered. The variable $x$ is proportional to the position along the channel, $t$ is the time variable, while $\eta(x,t), u(x,t)$ represent the deviation of the free surface from its undisturbed level at $(x,t)$ and the horizontal velocity of the fluid at $(x,t)$, respectively, see \cite{bonacs1,bonacs2} for details.

One of the relevant properties of (\ref{bse}) is the existence of solitary-wave solutions. They are solutions of traveling-wave form $u(X)=u(x-c_{s}t),\quad \eta(X)=\eta(x-c_{s}t)$ with constant speed $c_{s}>0$. The resulting waves are smooth, positive, even and decay to zero at infinity. By substitution into (\ref{bse}) and after one integration, the profiles $u=u(X), \eta=\eta(X)$ must satisfy the system
\begin{equation}\label{E422}
\begin{array}{l}
F\begin{pmatrix}
    u \cr
    \eta\end{pmatrix}=S\begin{pmatrix}
    u \\
    \eta\\
    \end{pmatrix}-\begin{pmatrix}
    u \eta \\
    \frac{u^{2}}{2}\\
    \end{pmatrix}
    =\begin{pmatrix}
    0 \cr
    0\end{pmatrix},\\
    S=\begin{pmatrix}
    -1&c_{s}(1-b\partial_{XX}) \\
    c_{s}(1-d\partial_{XX})&-(1+c\partial_{XX})\\
    \end{pmatrix}
\end{array}
\end{equation}
where $\partial_{XX}$ stands for the second order differentiation operator.
System (\ref{E422}) admits the group of spatial translations
\begin{equation}\label{E423}
G_{\alpha}(u,\eta)(X)=(u(X-\alpha),\eta(X-\alpha)),\quad \alpha\in\mathbb{R}
\end{equation}
as symmetry group. On the other hand, existence of these solutions is established from \cite{toland1}, (see also \cite{dougalism}) for $c_{s}>1$, although only in some cases analytical expression for the solitary waves is known \cite{chen}. Thus, for a typical value of $c_{s}$, the profiles $u,\eta$ must be approximated by numerical means.

In order to treat (\ref{E422}) numerically, considered here is the corresponding periodic problem on a sufficiently long interval $(-L,L)$. Then, (\ref{E422}) is discretized by using a Fourier collocation method, \cite{boyd}. This leads to the discrete system
\begin{equation}\label{E424}
\begin{array}{l}
F_{h}\begin{pmatrix}
    u_{h} \cr
    \eta_{h}\end{pmatrix}={S_{h}}\begin{pmatrix}
    u_{h} \\
    \eta_{h}\\
    \end{pmatrix}-\begin{pmatrix}
    u_{h}.\eta_{h} \\
    \frac{u_{h}.^{2}}{2}\\
    \end{pmatrix}=\begin{pmatrix}
    0 \cr
    0\end{pmatrix}\\
S_{h}=
\begin{pmatrix}
    -I_{N}&c_{s}(1-bD_{N}^{2})  \\
   c_{s}(1-dD_{N}^{2})&-(1+cD_{N}^{2})\\
    \end{pmatrix}
\end{array}
\end{equation}
where:
\begin{itemize}
\item $I_{N}$ stands for the $N\times N$ identity matrix and $D_{N}$ is the pseudospectral differentiation matrix, \cite{canutohqz}.
\item $\eta_{h}=\{\eta_{h,j}\}_{j}, u_{h}=\{u_{h,j}\}_{j}\in \mathbb{R}^{N}$ are vector approximations to $\widetilde{\eta}=\eta(x_{j}),
\widetilde{u}=u(x_{j})$ at the collocation points $x_{j}=-L+jh, h=2L/N, j=0,\ldots N-1$.
\end{itemize}
Note that (\ref{E424}) inherits a symmetry group, infinitesimally generated by the vector field $v:\mathbb{R}^{2N}\rightarrow\mathbb{R}^{2N}$,
\begin{eqnarray*}
(u_{h},\eta_{h})^{T}\mapsto v(u_{h},\eta_{h})=(D_{N}\eta_{h},D_{N}u_{h})^{T},
\end{eqnarray*}
which can be considered as a discrete version of (\ref{E423}), see the examples below.
For fixed $c_{s}$, an iterative technique can be applied to (\ref{E424}) in order to obtain the approximations $u_{h}, \eta_{h}$. This has been carried out by using the Newton's method
 \begin{eqnarray}
F_{h}^{\prime}(x_{n})\Delta(x_{n})&=&-F_{h}(x_{n}),\label{newt1}\\
x_{n+1}&=&x_{n}+\Delta x_{n},\quad n=0,1,\ldots,\label{newt2}
\end{eqnarray}
(where $x_{n}=(u_{n},\eta_{n})$). (The previously described Petviashvili method could also being applied.) Note that since $F_{h}^{\prime}(u_{h},\eta_{h})$ is singular, system (\ref{newt1}) may be ill-conditioned.
Two techniques for the resolution have been considered. The first one obtains $\Delta x_{n}=(\Delta u_{n},\Delta \eta_{n})$ that minimizes the Euclidean norm of the residual $r_{n}=-F_{h}(x_{n})-F_{h}^{\prime}(x_{n})\Delta x_{n}$ by using the minimum residual algorithm MINRES, \cite{paiges}. The method is implemented in physical space, by computing $D_{N}^{2}$ in (\ref{E424}) and the Jacobian
\begin{equation}\label{E425}
F_{h}^{\prime}(u_{n},\eta_{n})=S_{h}-\begin{pmatrix} \eta_{n}&u_{n}\\u_{n}&0\end{pmatrix},
\end{equation}
at each iteration. (Note that the Jacobian is symmetric.)  The second technique also makes use of the symmetry of (\ref{E425}) and solves (\ref{newt1}) with preconditioned conjugate gradient (PCG) techniques, \cite{GolubV,Meyer}. Although the PCG method can be initially applied to the symmetric, positive definite case, it also may be implemented when the matrix is symmetric but indefinite, \cite{demmel,yang1}.
In this case the implementation of the Newton's method is carried out in Fourier space and involves two iteration loops: an inner iteration corresponding to the resolution of (\ref{newt1}) with (PCG) and the outer iteration  associated to the advance (\ref{newt2}), \cite{lakoba,lakobay,yang1}. On the other hand, a discrete version of preconditioning operators of the form
\begin{eqnarray}
M=sI-\partial_{XX},\label{PO}
\end{eqnarray}
is used, where $s$ is a positive parameter to be adjusted according to the problem under implementation. The operators (\ref{PO}) are usually considered for this kind of problems, \cite{yang1,yang2}. The resulting method, consisting of the Newton's method (\ref{newt1}), (\ref{newt2}) along with PCG for (\ref{newt1}) will be denoted by (PCGN).

Both techniques have been performed and compared in terms of the behaviour of the residual error
\begin{eqnarray}\label{E4211}
RE_{n}=||F_{h}(u_{n},\eta_{n})||,\quad n=0,1,\ldots,
\end{eqnarray}
in the Euclidean norm, when $F_{h}$ is given by (\ref{E424}). This error and the corresponding iteration are controlled by a fixed tolerance $TOL=1E-14$.

Our experiments show that both methods are comparable in efficiency with respect to the number of iterations, but the second one (PCGN) provides a higher performance in terms of the computational time. Then the results presented here correspond to this second procedure.

We first compare the numerical results with an exact profile. As mentioned before, analytical expressions for solitary wave solutions of (\ref{E422}) and fixed $c_{s}$ are in general unknown. However, some profiles can be found in closed form. For example, \cite{chen}, for each $\theta^{2}\in (7/9,1)$, there is one solitary wave solution $(u_{s},\eta_{s})$ of the form
\begin{equation}\label{E426}
\begin{array}{l}
\eta_{s}(x,t)=\eta_{0}{\rm sech}^{2}(\lambda(x-c_{s}t-x_{0})),\\
u_{s}(x,t)=B\eta_{s}(x,t),
\end{array}
\end{equation}
where $x_{0}\in\mathbb{R}$ is arbitrary (it is the translational parameter; associated to the symmetry group) and the rest of the parameters is given by
\begin{eqnarray*}
&&\eta_{0}=\frac{9}{2}\frac{\theta^{2}-7/9}{1-\theta^{2}},\quad c_{s}=\frac{4(\theta^{2}-2/3)}{\sqrt{2(1-\theta^{2})(\theta^{2}-1/3)}},\\
&&\lambda=\frac{1}{2}\sqrt{\frac{3(\theta^{2}-7/9)}{(\theta^{2}-1/3)(\theta^{2}-2/3)}},\quad B=\sqrt{\frac{2(1-\theta^{2})}{\theta^{2}-1/3}}.
\end{eqnarray*}
\begin{figure}[htbp]
\centering \subfigure[]{
\includegraphics[width=0.45\textwidth]{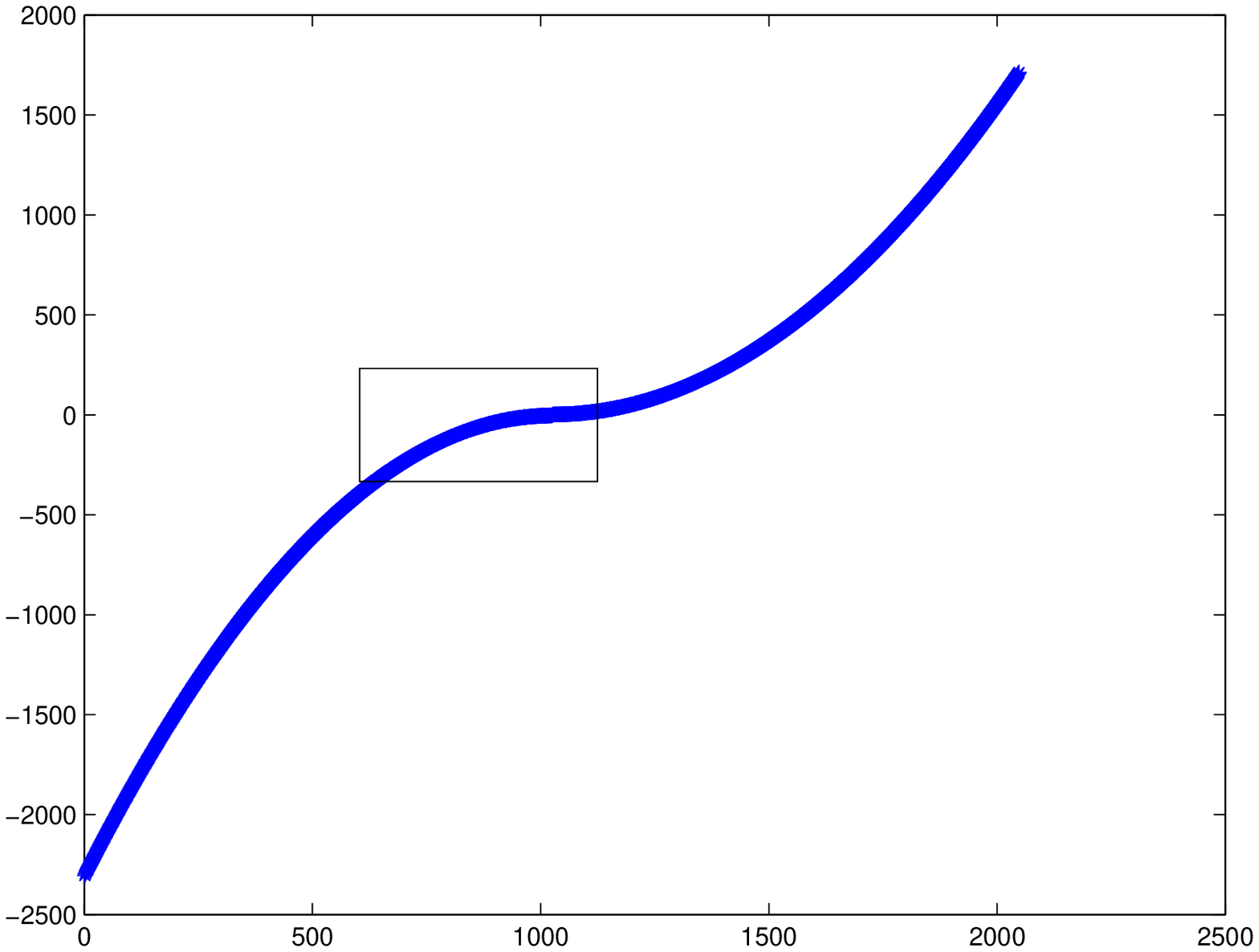} }
\subfigure[]{
\includegraphics[width=0.45\textwidth]{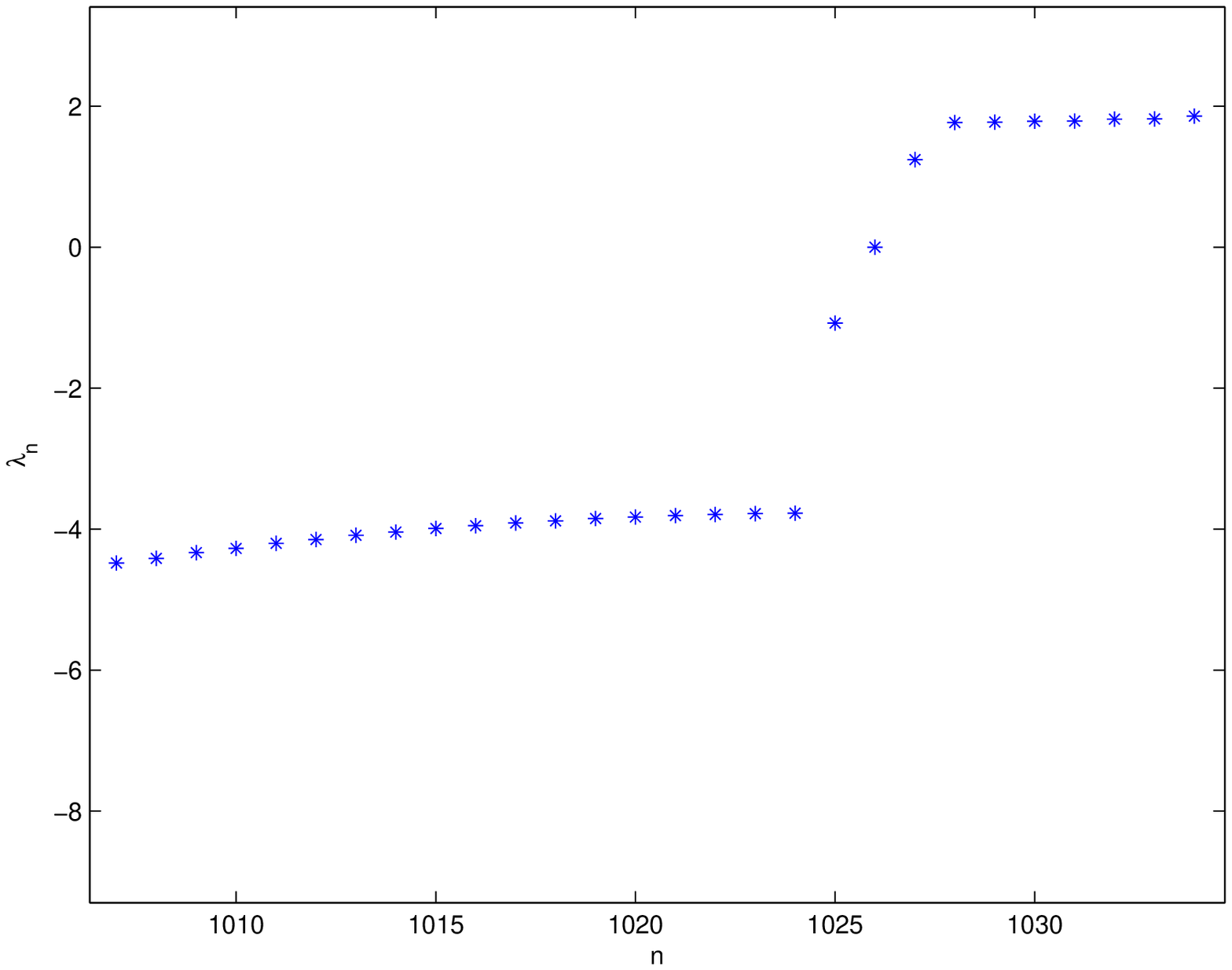} }
\caption{(a) Eigenvalues of the Jacobian (\ref{E425}) evaluated at the vector of exact profiles (\ref{E426}) at the grid points, $\theta^{2}=0.9, N=1024$. (b) Magnification of (a), with the eigenvalues between numbers $1007$ and $1034$, with $\mu=0$ at the position $1026$.} \label{examp2f1}
\end{figure}

Figure \ref{examp2f1}(a) shows the spectrum of the Jacobian (\ref{E425}) at the exact profiles (\ref{E426}) with $\theta^{2}=0.9$ (evaluated at the $N=1024$ Fourier collocation points; the size of the matrix is $2N\times 2N$). The indefinite character is observed and the magnification given by Figure \ref{examp2f1}(b) reveals the simple eigenvalue $\lambda=0$, associated to the translational symmetry of (\ref{E424}).
\begin{figure}[htbp]
\centering \subfigure[]{
\includegraphics[width=10cm]{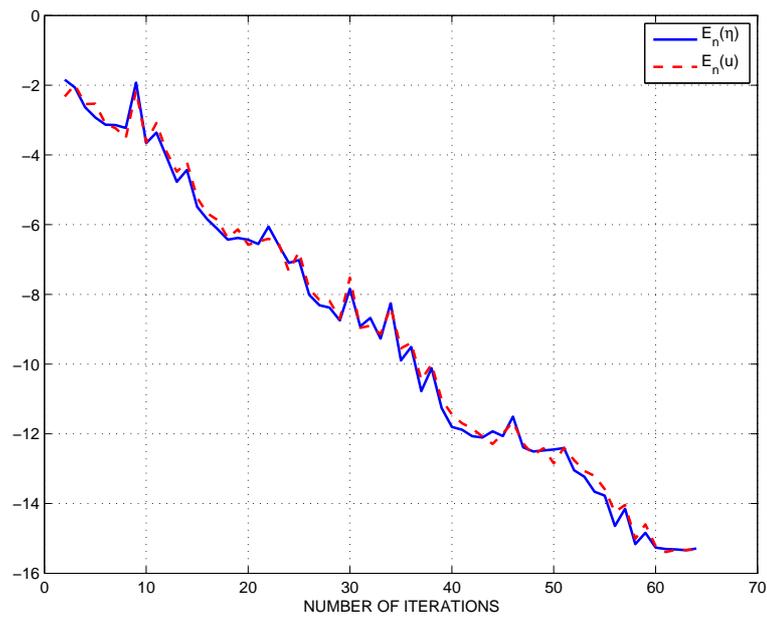} }
\subfigure[]{
\includegraphics[width=10cm]{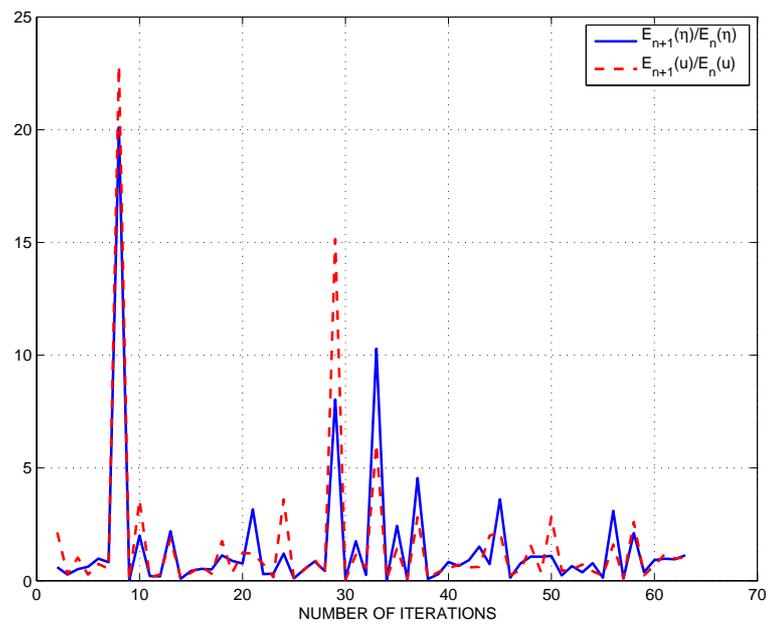} }
\caption{(a) Errors (\ref{E427}) for $\theta^{2}=0.9$ vs number of iterations, semi-logarithmic scale. (b) Ratios $E_{n+1}(\eta)/E_{n}(\eta), E_{n+1}(u)/E_{n}(u)$ vs number of iterations.} \label{examp2f2}
\end{figure}

The results of the method are first compared with the exact profiles in Figure \ref{examp2f2}(a). This shows the behaviour of the error between the iteration and the profiles
\begin{eqnarray}
\label{E427}
E_{n}(\eta)=||\eta_{n}-\widetilde{\eta}||,\quad E_{n}(u)=||u_{n}-\widetilde{u}||,\quad n=0,1,\ldots
\end{eqnarray}
in logarithmic scale and as function of the number of iterations. (The Euclidean norm is considered.) The initial iteration is a perturbation of the exact profiles of the form
\begin{equation}\label{E428}
\begin{array}{l}
\eta_{0}(X)=\eta_{s}(X)+\epsilon e^{-X^{2}},\\
u_{0}(X)=u_{s}(X)+\epsilon e^{-X^{2}},
\end{array}
\end{equation}
with $\epsilon=5E-02$. The local convergence of the method is illustrated. On the other hand, the order of this convergence is displayed in Figure \ref{examp2f2}(b). This shows the ratios $E_{n+1}(\eta)/E_{n}(\eta), E_{n+1}(u)/E_{n}(u)$, also as function of the number $n$ of iterations. These ratios are bounded with $n$, suggesting a linear convergenceand they look to rule out the $Q$ convergence, \cite{dennism1,dennism2}.
\begin{figure}[htbp]
\centering \subfigure[]{
\includegraphics[width=10cm]{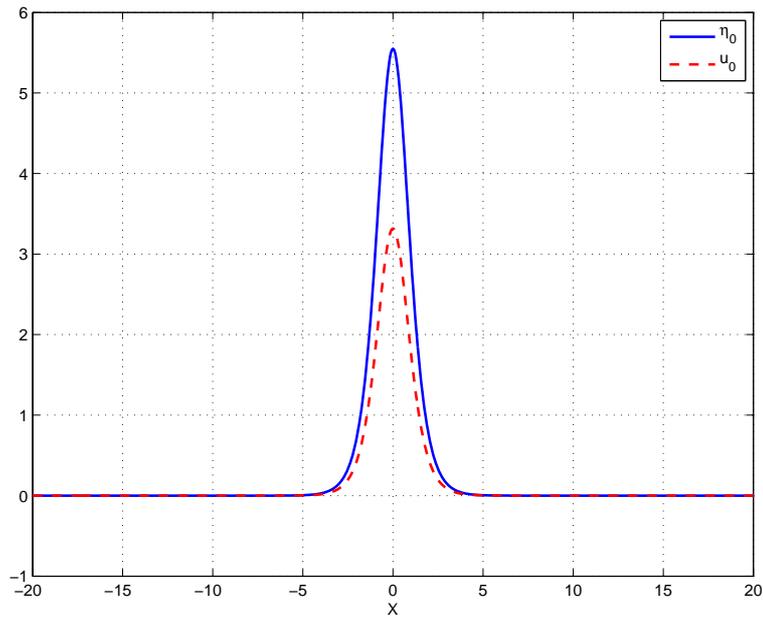} }
\subfigure[]{
\includegraphics[width=10cm]{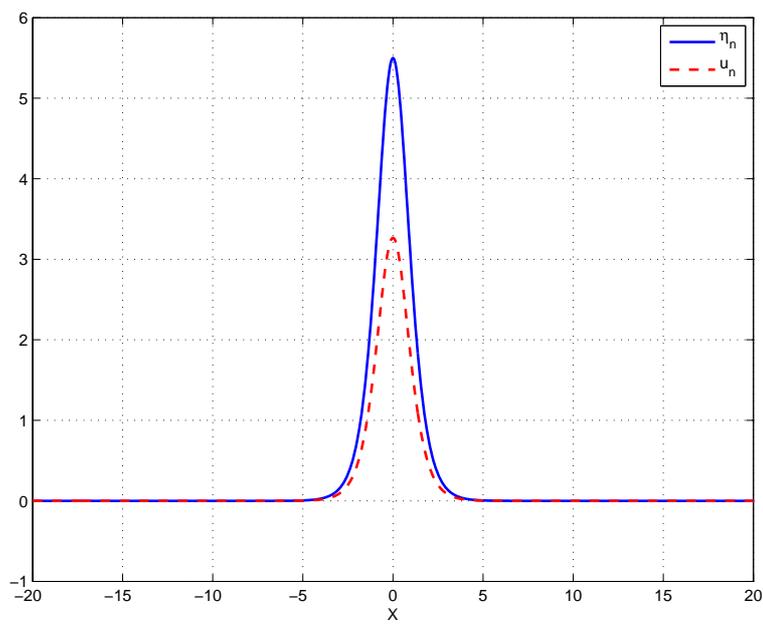} }
\caption{(a) Initial iteration of the form (\ref{E428}) with $\epsilon=5E-02$. (b) Final iterate obtained by the (PCGN) method.} \label{examp2f3}
\end{figure}

The orbital convergence is illustrated by the following figures. Figure \ref{examp2f3}(b) shows the final approximate profiles obtained by the (PCGN) method, with initial iteration displayed in Figure \ref{examp2f3}(a). This is of the form (\ref{E428}) with $\epsilon=5E-02$. The initial profiles are centered at the origin, and the displacement from it of the final iterates are computed as $x_{\eta}=1.235296E-15, x_{u}=2.334808E-15$, so they can be considered as centered at the origin as well.
\begin{figure}[htbp]
\centering \subfigure[]{
\includegraphics[width=10cm]{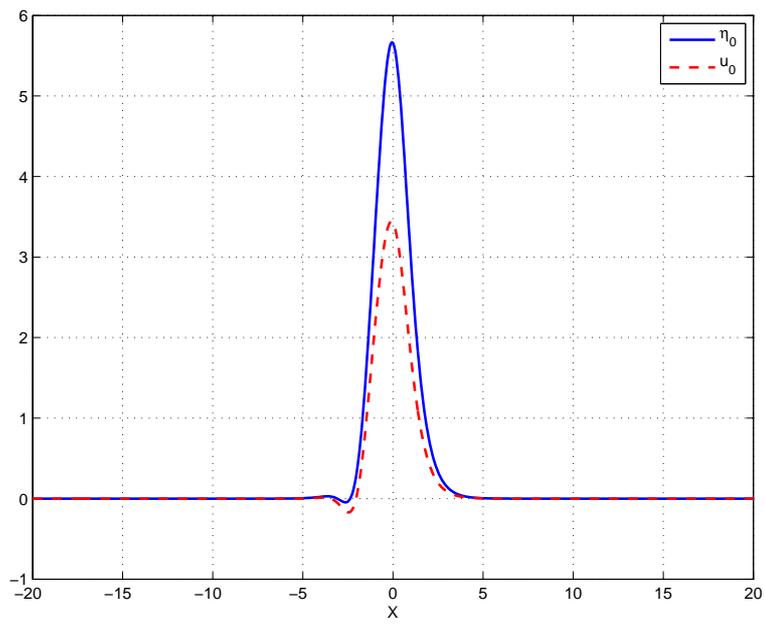} }
\subfigure[]{
\includegraphics[width=10cm]{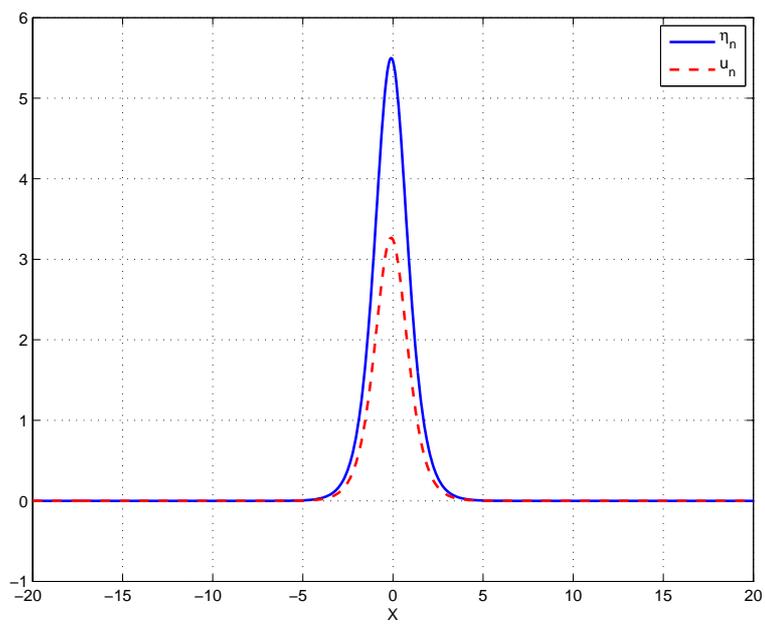} }
\caption{(a) Initial iteration of the form (\ref{E429}) with $X_{0}=-1.5, \epsilon=1$. (b) Final iterate obtained by the (PCGN) method.} \label{examp2f4}
\end{figure}
The effect of the translational symmetry is observed in Figure  \ref{examp2f4}. From a profile of the form
\begin{equation}\label{E429}
\begin{array}{l}
\eta_{0}(X)=\eta_{s}(X)+\epsilon (X-X_{0})e^{-(X-X_{0})^{2}},\\
u_{0}(X)=u_{s}(X)+\epsilon (X-X_{0})e^{-(X-X_{0})^{2}},
\end{array}
\end{equation}
with $X_{0}=-1.5$ and $\epsilon=1$, shown in Figure  \ref{examp2f4}(a), the procedure converges (up to the fixed tolerance $1E-14$) with the final iterates shown in Figure  \ref{examp2f4}(b). The profiles are affected by a phase shift and centered at
\begin{eqnarray*}
x_{\eta}\sim x_{u}=-1.012779E-01.
\end{eqnarray*}
(The computed $x_{\eta}, x_{u}$ coincide up to ten digits.)

\begin{table}
\begin{center}
\begin{tabular}{|c|c|c|}
\hline\hline
$\epsilon$ & $x_{u}$&$x_{\eta}$\\\hline
1E-01&-9.9534E-02&-9.9534E-02\\
5E-02&-4.9941E-02&-4.9941E-02\\
1E-02&-9.9995E-03&-9.9995E-03\\
5E-03&-4.9999E-03&-4.9999E-03\\
\hline\hline
\end{tabular}
\end{center}
\caption{Phase shift of the final iterates given by the PCGN method from initial iterations (\ref{E4210}) and the values of $\epsilon$ displayed in the first column. For each row, $x_{\eta}, x_{u}$ coincide up to ten digits.}\label{tab_etna3}
\end{table}
A final experiment, related to the profiles (\ref{E426}), checks the local convergence and the effect of the components of the initial iteration associated to the symmetry on the limit profile. Initial iterations of the form
\begin{equation}\label{E4210}
\begin{array}{l}
\eta_{0}(X)=\eta_{s}(X)+\epsilon D_{n}\eta_{s},\\
u_{0}(X)=u_{s}(X)+\epsilon D_{N}u_{s},
\end{array}
\end{equation}
and several (small) values of $\epsilon$, the phase shift of the final iterates have been computed and are shown in Table \ref{tab_etna3}. The resulting displacements $x_{\eta}, x_{u}$ correspond, as $\epsilon$ decreases, to the limit profile $G_{\epsilon}(u_{h},\eta_{h})$ associated to a spatial translation of $\epsilon$.
\begin{figure}[htbp]
\centering \subfigure[]{
\includegraphics[width=10cm]{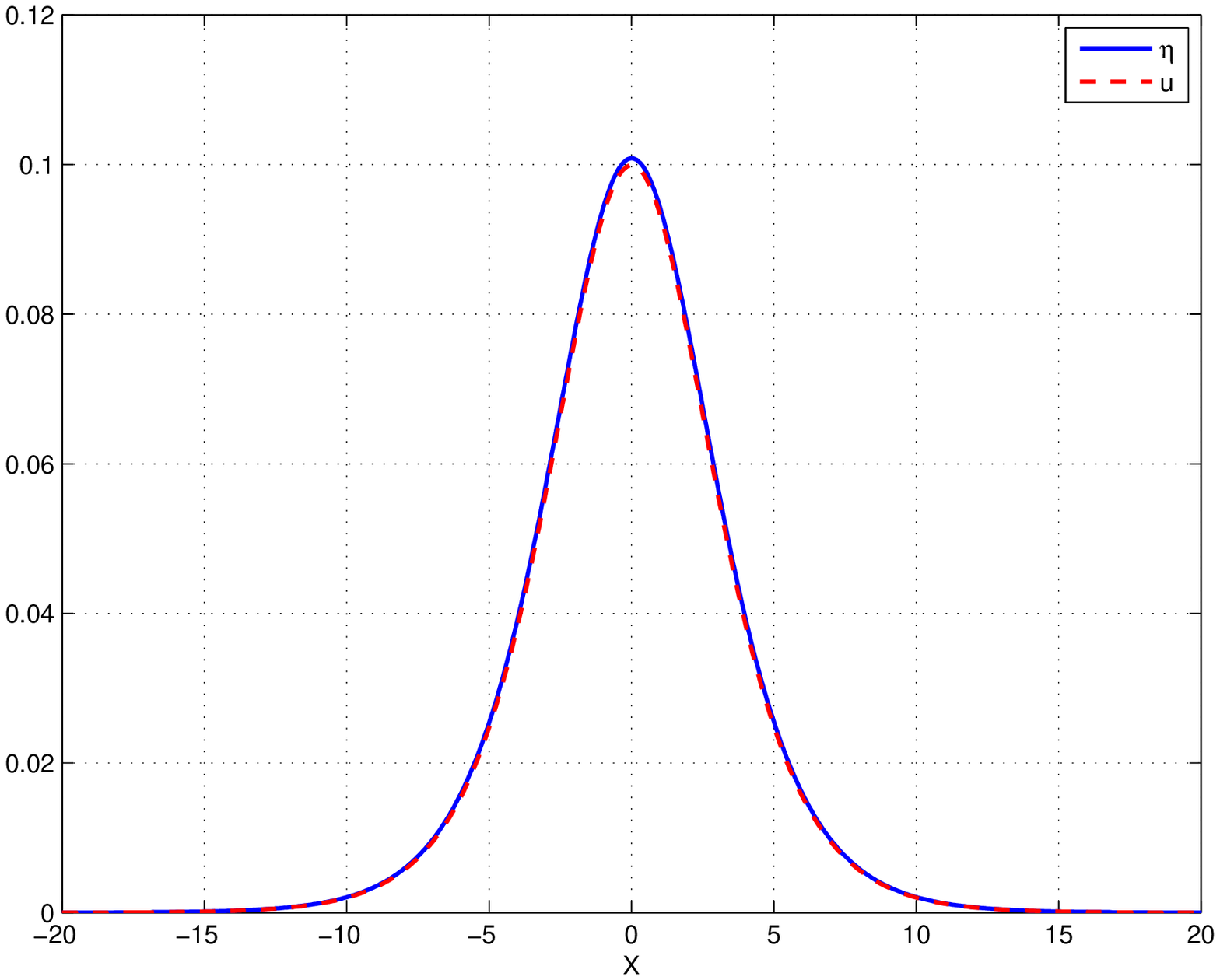} }
\subfigure[]{
\includegraphics[width=10cm]{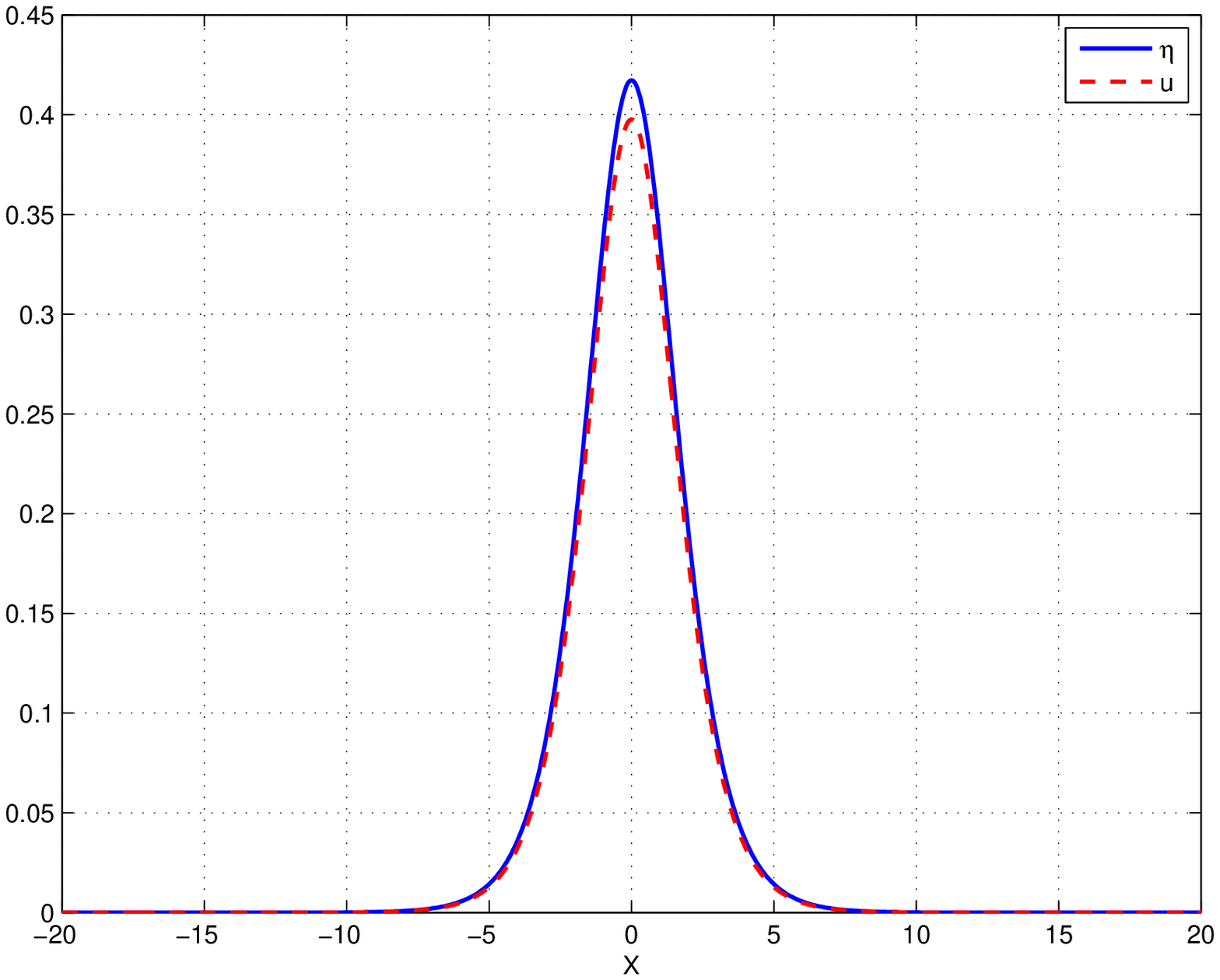} }
\caption{Final iterates with the (PCGN) method and from initial profiles (\ref{E426}) with $\theta^{2}=0.9$ for (a) $c_{s}=1.05$; (b) $c_{s}=1.2$.} \label{examp2f6}
\end{figure}

A second group of experiments checks the accuracy of the method to approximate solitary wave profile solutions of (\ref{E422}) for which explicit formulas are not known. We have taken $\theta^{2}=0.9$ and two speeds, $c_{s}=1.05, 1.2$. By using as initial iterations the profiles (\ref{E426}) with $\theta^{2}=0.9$, the final iterates  are shown in Figure \ref{examp2f6}. For the first speed, Figure \ref{examp2f6}(a), the profiles are displaced with a phase shift (the same, up to ten digits, for both $\eta$ and $u$) of $x_{\eta}\sim x_{u}= 5.177124E-03$, while for the waves of Figure \ref{examp2f6}(b), corresponding to the second value of the speed, the computed phase shift, for both components, is below the machine precision, therefore no displacement is assumed.

\begin{figure}[htbp]
\centering
\includegraphics[width=10cm]{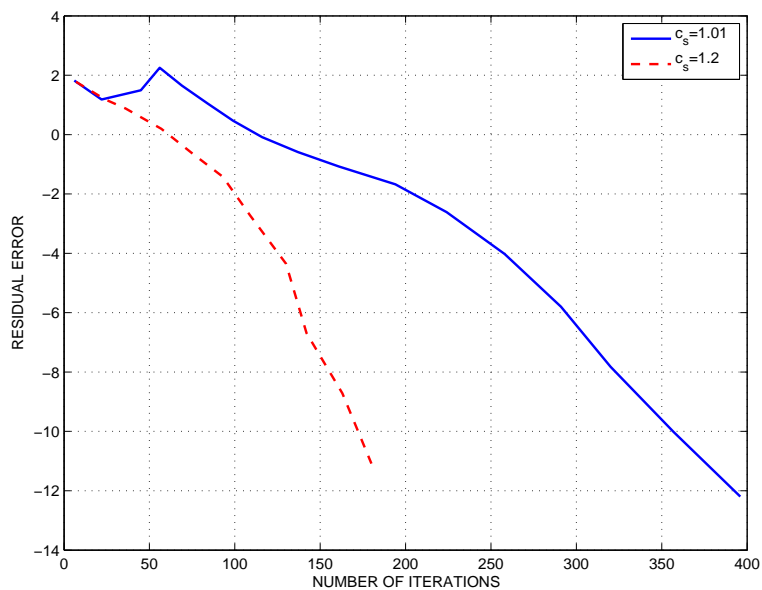}
\caption{Residual error (\ref{E4211}) for the computed profiles shown in Figure \ref{examp2f6}.} \label{examp2f7}
\end{figure}
\begin{figure}[htbp]
\centering
\subfigure[]{
\includegraphics[width=10cm]{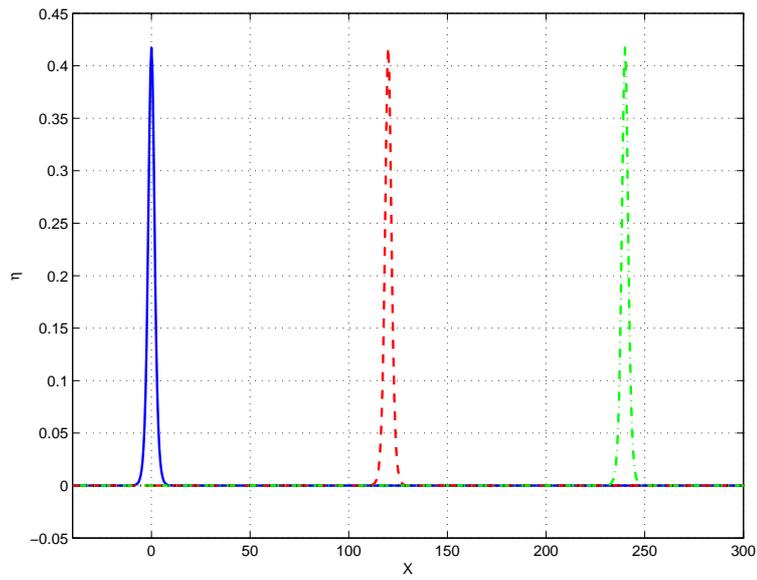} }
\subfigure[]{
\includegraphics[width=10cm]{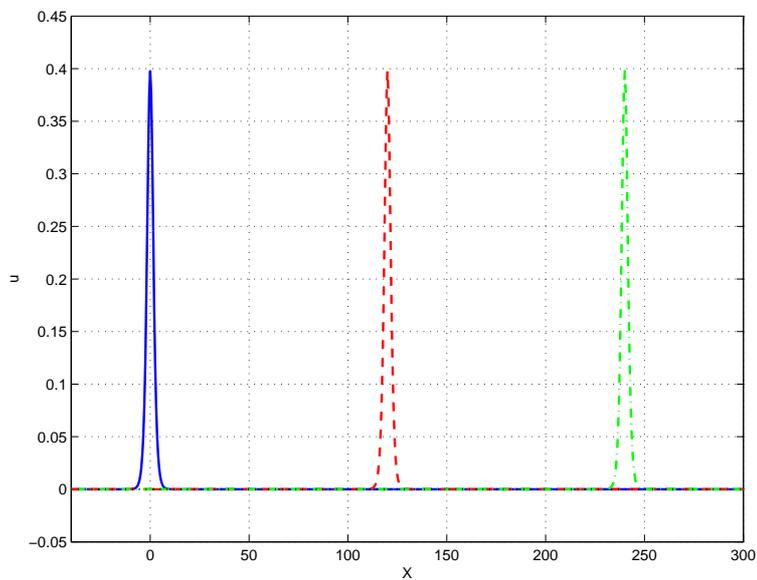} }
\caption{Form of the numerical solution of a time-stepping code for (\ref{bse}) with initial condition given by the profiles of Figure \ref{examp2f6}(b) at times $t=0, 100, 200$ (solid, dashed and dashed-dotted lines, respectively). (a) $\eta$ profile; (b) $u$ profile.} \label{examp2f8}
\end{figure}

The accuracy of the computed waves is checked in Figures \ref{examp2f7} and \ref{examp2f8}. Figure \ref{examp2f7} shows the residual error (\ref{E4211}) as function of the number of iterations and in logarithmic scale, for the computed profiles displayed in Figure \ref{examp2f6}. As is known, (\ref{E4211}) measures the error generated by the iterations in the formula (\ref{E424}) for the profiles. The method gets a residual error of order of about $10^{-12}$ in less iterations in the case $c_{s}=1.2$ than for $c_{s}=1.05$. (The number of iterations depends, among others, on the preconditioning parameter $s$ in (\ref{PO}). An optimal value of it may be different for different speeds.)

Finally, a necessary test to check the accuracy of the computed profiles as solitary wave solutions of the Bona-Smith system (\ref{bse}) consists of taking them as initial conditions of a code to integrate (\ref{bse}) numerically. Figure \ref{examp2f8} illustrates the evolution of the solution, by displaying the form of the corresponding approximation at different times. (The experiment corresponds to $c_{s}=1.2$.) A propagation without backward or forward disturbances is observed, a solitary wave-like evolution to the right with speed $c_{s}=1.2$.

\end{document}